\documentclass[12pt]{article}
\usepackage{latexsym,amssymb,amsmath,ulem}
\usepackage{color}
\usepackage{graphicx}
\usepackage{latexsym,amssymb,amsmath,ulem}
\usepackage{amsfonts}
\usepackage{latexsym}
\usepackage{graphicx}
\usepackage{color}

\def\tto{\;{\lower 1pt \hbox{$\rightarrow$}}\kern -10pt
\hbox{\raise 2pt \hbox{$\rightarrow$}}\;}

\newtheorem{theorem}{Theorem}[section]
\newtheorem{proposition}{Proposition}[section]
\newtheorem{corollary}{Corollary}[section]
\newtheorem{lemma}{Lemma}[section]
\newtheorem{remark}{Remark}[section]
\newtheorem{example}{Example}[section]
\newtheorem{definition}{Definition}[section]

\numberwithin{equation}{section}

\renewcommand{\theequation}{\thesection.\arabic{equation}}
\normalsize

\makeatletter
\providecommand{\BOXEDSPECIAL}[4]{\hbox to #2{\raise #3\hbox to #2{\null #1\hfil}}}

\newcount\@hour\newcount\@minute\chardef\@x10\chardef\@xv60
\def\tcitime{
\def\@time{%
  \@minute\time\@hour\@minute\divide\@hour\@xv
  \ifnum\@hour<\@x 0\fi\the\@hour:%
  \multiply\@hour\@xv\advance\@minute-\@hour
  \ifnum\@minute<\@x 0\fi\the\@minute
  }}%

\@ifundefined{hyperref}{}{}

\@ifundefined{qExtProgCall}{\def\qExtProgCall#1#2#3#4#5#6{\relax}}{}
\def\QCTOpt[#1]#2{%
  \def\QCTOptB{#1}
  \def\QCTOptA{#2}
}
\def\QCTNOpt#1{%
  \def\QCTOptA{#1}
  \let\QCTOptB\empty
}
\def\Qct{%
  \@ifnextchar[{%
    \QCTOpt}{\QCTNOpt}
}
\def\QCBOpt[#1]#2{%
  \def\QCBOptB{#1}
  \def\QCBOptA{#2}
}
\def\QCBNOpt#1{%
  \def\QCBOptA{#1}
  \let\QCBOptB\empty
}
\def\Qcb{%
  \@ifnextchar[{%
    \QCBOpt}{\QCBNOpt}
}
\def\PrepCapArgs{%
  \ifx\QCBOptA\empty
    \ifx\QCTOptA\empty
      {}%
    \else
      \ifx\QCTOptB\empty
        {\QCTOptA}%
      \else
        [\QCTOptB]{\QCTOptA}%
      \fi
    \fi
  \else
    \ifx\QCBOptA\empty
      {}%
    \else
      \ifx\QCBOptB\empty
        {\QCBOptA}%
      \else
        [\QCBOptB]{\QCBOptA}%
      \fi
    \fi
  \fi
}
\newcount\GRAPHICSTYPE
\GRAPHICSTYPE=\z@
\def\GRAPHICSPS#1{%
 \ifcase\GRAPHICSTYPE
   \special{ps: #1}%
 \or
   \special{language "PS", include "#1"}%
 \fi
}%
\def\graffile#1#2#3#4#5{%
    \bgroup
    \leavevmode
    \@ifundefined{bbl@deactivate}{\def~{\string~}}{\activesoff}
    \raise -#4 \BOXTHEFRAME{%
       \BOXEDSPECIAL{#1}{#2}{#3}{#5}}%
    \egroup
}%
\def\draftbox#1#2#3#4{%
 \leavevmode\raise -#4 \hbox{%
  \frame{\rlap{\protect\tiny #1}\hbox to #2%
   {\vrule height#3 width\z@ depth\z@\hfil}%
  }%
 }%
}%
\newcount\draft
\draft=\z@

\newif\ifwasdraft
\wasdraftfalse

\def\GRAPHIC#1#2#3#4#5{%
 \ifnum\draft=\@ne\draftbox{#2}{#3}{#4}{#5}%
  \else\graffile{#1}{#3}{#4}{#5}{#2}%
  \fi
 }%
\def\addtoLaTeXparams#1{%
    \edef\LaTeXparams{\LaTeXparams #1}}%
\newif\ifBoxFrame \BoxFramefalse
\newif\ifOverFrame \OverFramefalse
\newif\ifUnderFrame \UnderFramefalse

\def\BOXTHEFRAME#1{%
   \hbox{%
      \ifBoxFrame
         \frame{#1}%
      \else
         {#1}%
      \fi
   }%
}

\def\doFRAMEparams#1{\BoxFramefalse\OverFramefalse\UnderFramefalse\readFRAMEparams#1\end}%
\def\readFRAMEparams#1{%
 \ifx#1\end%
  \let\next=\relax
  \else
  \ifx#1i\dispkind=\z@\fi
  \ifx#1d\dispkind=\@ne\fi
  \ifx#1f\dispkind=\tw@\fi
  \ifx#1t\addtoLaTeXparams{t}\fi
  \ifx#1b\addtoLaTeXparams{b}\fi
  \ifx#1p\addtoLaTeXparams{p}\fi
  \ifx#1h\addtoLaTeXparams{h}\fi
  \ifx#1X\BoxFrametrue\fi
  \ifx#1O\OverFrametrue\fi
  \ifx#1U\UnderFrametrue\fi
  \ifx#1w
    \ifnum\draft=1\wasdrafttrue\else\wasdraftfalse\fi
    \draft=\@ne
  \fi
  \let\next=\readFRAMEparams
  \fi
 \next
 }%
\def\IFRAME#1#2#3#4#5#6{%
      \bgroup
      \let\QCTOptA\empty
      \let\QCTOptB\empty
      \let\QCBOptA\empty
      \let\QCBOptB\empty
      #6%
      \parindent=0pt%
      \leftskip=0pt
      \rightskip=0pt
      \setbox0 = \hbox{\QCBOptA}%
      \@tempdima = #1\relax
      \ifOverFrame
          \typeout{This is not implemented yet}%
          \show\HELP
      \else
         \ifdim\wd0>\@tempdima
            \advance\@tempdima by \@tempdima
            \ifdim\wd0 >\@tempdima
               \textwidth=\@tempdima
               \setbox1 =\vbox{%
                  \noindent\hbox to \@tempdima{\hfill\GRAPHIC{#5}{#4}{#1}{#2}{#3}\hfill}\\%
                  \noindent\hbox to \@tempdima{\parbox[b]{\@tempdima}{\QCBOptA}}%
               }%
               \wd1=\@tempdima
            \else
               \textwidth=\wd0
               \setbox1 =\vbox{%
                 \noindent\hbox to \wd0{\hfill\GRAPHIC{#5}{#4}{#1}{#2}{#3}\hfill}\\%
                 \noindent\hbox{\QCBOptA}%
               }%
               \wd1=\wd0
            \fi
         \else
            \ifdim\wd0>0pt
              \hsize=\@tempdima
              \setbox1 =\vbox{%
                \unskip\GRAPHIC{#5}{#4}{#1}{#2}{0pt}%
                \break
                \unskip\hbox to \@tempdima{\hfill \QCBOptA\hfill}%
              }%
              \wd1=\@tempdima
           \else
              \hsize=\@tempdima
              \setbox1 =\vbox{%
                \unskip\GRAPHIC{#5}{#4}{#1}{#2}{0pt}%
              }%
              \wd1=\@tempdima
           \fi
         \fi
         \@tempdimb=\ht1
         \advance\@tempdimb by \dp1
         \advance\@tempdimb by -#2%
         \advance\@tempdimb by #3%
         \leavevmode
         \raise -\@tempdimb \hbox{\box1}%
      \fi
      \egroup%
}%
\def\DFRAME#1#2#3#4#5{%
 \begin{center}
     \let\QCTOptA\empty
     \let\QCTOptB\empty
     \let\QCBOptA\empty
     \let\QCBOptB\empty
     \ifOverFrame
        #5\QCTOptA\par
     \fi
     \GRAPHIC{#4}{#3}{#1}{#2}{\z@}
     \ifUnderFrame
        \nobreak\par\nobreak#5\QCBOptA
     \fi
 \end{center}%
 }%
\def\FFRAME#1#2#3#4#5#6#7{%
 \begin{figure}[#1]%
  \let\QCTOptA\empty
  \let\QCTOptB\empty
  \let\QCBOptA\empty
  \let\QCBOptB\empty
  \ifOverFrame
    #4
    \ifx\QCTOptA\empty
    \else
      \ifx\QCTOptB\empty
        \caption{\QCTOptA}%
      \else
        \caption[\QCTOptB]{\QCTOptA}%
      \fi
    \fi
    \ifUnderFrame\else
      \label{#5}%
    \fi
  \else
    \UnderFrametrue%
  \fi
  \begin{center}\GRAPHIC{#7}{#6}{#2}{#3}{\z@}\end{center}%
  \ifUnderFrame
    #4
    \ifx\QCBOptA\empty
      \caption{}%
    \else
      \ifx\QCBOptB\empty
        \caption{\QCBOptA}%
      \else
        \caption[\QCBOptB]{\QCBOptA}%
      \fi
    \fi
    \label{#5}%
  \fi
  \end{figure}%
 }%
\newcount\dispkind%

\def\makeactives{
  \catcode`\"=\active
  \catcode`\;=\active
  \catcode`\:=\active
  \catcode`\'=\active
  \catcode`\~=\active
}
\bgroup
   \makeactives
   \gdef\activesoff{%
      \def"{\string"}
      \def;{\string;}
      \def:{\string:}
      \def'{\string'}
    }
\egroup

\def\FRAME#1#2#3#4#5#6#7#8{%
 \bgroup
 \ifnum\draft=\@ne
   \wasdrafttrue
 \else
   \wasdraftfalse%
 \fi
 \def\LaTeXparams{}%
 \dispkind=\z@
 \def\LaTeXparams{}%
 \doFRAMEparams{#1}%
 \ifnum\dispkind=\z@\IFRAME{#2}{#3}{#4}{#7}{#8}{#5}\else
  \ifnum\dispkind=\@ne\DFRAME{#2}{#3}{#7}{#8}{#5}\else
   \ifnum\dispkind=\tw@
    \edef\@tempa{\noexpand\FFRAME{\LaTeXparams}}%
    \@tempa{#2}{#3}{#5}{#6}{#7}{#8}%
    \fi
   \fi
  \fi
  \ifwasdraft\draft=1\else\draft=0\fi{}%
  \egroup
 }%
\def\TEXUX#1{"texux"}

\def\limfunc#1{\mathop{\rm #1}}%
\def\func#1{\mathop{\rm #1}\nolimits}%
\long\def\QQQ#1#2{%
     \long\expandafter\def\csname#1\endcsname{#2}}%
\@ifundefined{QTP}{\def\QTP#1{}}{}
\@ifundefined{QEXCLUDE}{\def\QEXCLUDE#1{}}{}
\@ifundefined{Qlb}{}{}
\@ifundefined{Qlt}{}{}
\long\def\QQA#1#2{}%
\def\EXPAND#1[#2]#3{}%
\def\NOEXPAND#1[#2]#3{}%
\def\LaTeXparent#1{}%
\def\ChildStyles#1{}%
\def\ChildDefaults#1{}%
\def\QTagDef#1#2#3{}%

\@ifundefined{correctchoice}{}{}
\@ifundefined{HTML}{\def\HTML#1{\relax}}{}
\@ifundefined{TCIIcon}{\def\TCIIcon#1#2#3#4{\relax}}{}
\if@compatibility
\else
  \providecommand{\UNICODE}[2][]{}
  
\fi

\@ifundefined{StyleEditBeginDoc}{}{}
\def\QQfnmark#1{\footnotemark}

\@ifundefined{TCIMAKEINDEX}{}{\makeindex}%
\@ifundefined{abstract}{%
 \def\abstract{%
  \if@twocolumn
   \section*{Abstract (Not appropriate in this style!)}%
   \else \small
   \begin{center}{\bf Abstract\vspace{-.5em}\vspace{\z@}}\end{center}%
   \quotation
   \fi
  }%
 }{%
 }%
\@ifundefined{endabstract}{\def\endabstract
  {\if@twocolumn\else\endquotation\fi}}{}%
\@ifundefined{maketitle}{\def\maketitle#1{}}{}%
\@ifundefined{affiliation}{\def\affiliation#1{}}{}%
\@ifundefined{proof}{}{}%
\@ifundefined{endproof}{}{}%
\@ifundefined{newfield}{\def\newfield#1#2{}}{}%
\@ifundefined{chapter}{\def\chapter#1{\par(Chapter head:)#1\par }%
 \newcount\c@chapter}{}%
\@ifundefined{part}{\def\part#1{\par(Part head:)#1\par }}{}%
\@ifundefined{section}{\def\section#1{\par(Section head:)#1\par }}{}%
\@ifundefined{subsection}{\def\subsection#1%
 {\par(Subsection head:)#1\par }}{}%
\@ifundefined{subsubsection}{\def\subsubsection#1%
 {\par(Subsubsection head:)#1\par }}{}%
\@ifundefined{paragraph}{\def\paragraph#1%
 {\par(Subsubsubsection head:)#1\par }}{}%
\@ifundefined{subparagraph}{\def\subparagraph#1%
 {\par(Subsubsubsubsection head:)#1\par }}{}%
\@ifundefined{therefore}{}{}%
\@ifundefined{backepsilon}{}{}%
\@ifundefined{yen}{}{}%
\@ifundefined{registered}{%
   \def\registered{\relax\ifmmode{}\r@gistered
                    \else$\m@th\r@gistered$\fi}%
 \def\r@gistered{^{\ooalign
  {\hfil\raise.07ex\hbox{$\scriptstyle\rm\text{R}$}\hfil\crcr
  \mathhexbox20D}}}}{}%
\@ifundefined{Eth}{}{}%
\@ifundefined{eth}{}{}%
\@ifundefined{Thorn}{}{}%
\@ifundefined{thorn}{}{}%
\@ifundefined{degree}{}{}%
\newdimen\theight
\def\Column{%
 \vadjust{\setbox\z@=\hbox{\scriptsize\quad\quad tcol}%
  \theight=\ht\z@\advance\theight by \dp\z@\advance\theight by \lineskip
  \kern -\theight \vbox to \theight{%
   \rightline{\rlap{\box\z@}}%
   \vss
   }%
  }%
 }%
\def\qed{%
 \ifhmode\unskip\nobreak\fi\ifmmode\ifinner\else\hskip5\p@\fi\fi
 \hbox{\hskip5\p@\vrule width4\p@ height6\p@ depth1.5\p@\hskip\p@}%
 }%
\def\miss{\hbox{\vrule height2\p@ width 2\p@ depth\z@}}%
\def\tcol#1{{\baselineskip=6\p@ \vcenter{#1}} \Column}  %
\@ifundefined{note}{}{}%

\def\newfmtname{LaTeX2e}

\ifx\fmtname\newfmtname
  \DeclareOldFontCommand{\rm}{\normalfont\rmfamily}{\mathrm}
  \DeclareOldFontCommand{\sf}{\normalfont\sffamily}{\mathsf}
  \DeclareOldFontCommand{\tt}{\normalfont\ttfamily}{\mathtt}
  \DeclareOldFontCommand{\bf}{\normalfont\bfseries}{\mathbf}
  \DeclareOldFontCommand{\it}{\normalfont\itshape}{\mathit}
  \DeclareOldFontCommand{\sl}{\normalfont\slshape}{\@nomath\sl}
  \DeclareOldFontCommand{\sc}{\normalfont\scshape}{\@nomath\sc}
\fi

\@ifundefined{theorem}{\newtheorem{theorem}{Theorem}}{}
\@ifundefined{lemma}{\newtheorem{lemma}[theorem]{Lemma}}{}
\@ifundefined{corollary}{\newtheorem{corollary}[theorem]{Corollary}}{}
\@ifundefined{conjecture}{}{}
\@ifundefined{proposition}{\newtheorem{proposition}[theorem]{Proposition}}{}
\@ifundefined{axiom}{}{}
\@ifundefined{remark}{\newtheorem{remark}{Remark}}{}
\@ifundefined{example}{\newtheorem{example}{Example}}{}
\@ifundefined{exercise}{}{}
\@ifundefined{definition}{\newtheorem{definition}{Definition}}{}

\@ifundefined{mathletters}{%
  \newcounter{equationnumber}
  \def\mathletters{%
     \addtocounter{equation}{1}
     \edef\@currentlabel{\theequation}%
     \setcounter{equationnumber}{\c@equation}
     \setcounter{equation}{0}%
     \edef\theequation{\@currentlabel\noexpand\alph{equation}}%
  }
  
}{}

\@ifundefined{BibTeX}{%
    \def\BibTeX{{\rm B\kern-.05em{\sc i\kern-.025em b}\kern-.08em
                 T\kern-.1667em\lower.7ex\hbox{E}\kern-.125emX}}}{}%
\@ifundefined{AmS}%
    {\def\AmS{{\protect\usefont{OMS}{cmsy}{m}{n}%
                A\kern-.1667em\lower.5ex\hbox{M}\kern-.125emS}}}{}%
\@ifundefined{AmSTeX}{}{}%
\def\@@eqncr{\let\@tempa\relax
    \ifcase\@eqcnt \def\@tempa{& & &}\or \def\@tempa{& &}%
      \else \def\@tempa{&}\fi
     \@tempa
     \if@eqnsw
        \iftag@
           \@taggnum
        \else
           \@eqnnum\stepcounter{equation}%
        \fi
     \fi
     \global\tag@false
     \global\@eqnswtrue
     \global\@eqcnt\z@\cr}

\def\TCItag{\@ifnextchar*{\@TCItagstar}{\@TCItag}}
\def\@TCItag#1{%
    \global\tag@true
    \global\def\@taggnum{(#1)}}
\def\@TCItagstar*#1{%
    \global\tag@true
    \global\def\@taggnum{#1}}
\def\QATOP#1#2{{#1 \atop #2}}%
\def\dsum{\mathop{\displaystyle \sum }}%
\def\dbigcup{\mathop{\displaystyle \bigcup }}%
\makeatother

\input{tcilatex}
\begin{document}

\title{Duality for the robust sum of functions}
\author{N. Dinh\thanks{
International University, Vietnam National University - HCMC, Linh Trung
ward, Thu Duc district, Ho Chi Minh city, Vietnam (ndinh02@gmail.com). }, \
\ M.A. Goberna\thanks{%
Department of Mathematics, University of Alicante, Spain (mgoberna@ua.es){}}%
, \ \ M. Volle\thanks{%
Avignon University, LMA EA 2151, Avignon, France
(michel.volle@univ-avignon.fr)} }
\maketitle
\date{}

\begin{abstract}
In this paper we associate with an infinite family of real extended
functions defined on a locally convex space, a sum, called robust sum, which
is always well-defined. We also associate with that family of functions a
dual pair of problems formed by the unconstrained minimization of its robust
sum and the so-called optimistic dual. For such a dual pair, we characterize
weak duality, zero duality gap, and strong duality, and their corresponding
stable versions, in terms of multifunctions associated with the given family
of functions and a given approximation parameter $\varepsilon \geq 0$ which
is related to the $\varepsilon $-subdifferential of the robust sum of the
family. We also consider the particular case when all functions of the
family are convex, assumption allowing to characterize the duality
properties in terms of closedness conditions.
\end{abstract}


\qquad \textbf{Keywords} \ Robust sum function \textperiodcentered\ Weak
duality \textperiodcentered\ Zero uality gap \textperiodcentered\ Strong

\qquad duality \textperiodcentered\ Stable duality theorems

\qquad \textbf{Mathematics Subject Classifications (2010)} \
90C46\textperiodcentered\ 49N15 \textperiodcentered\ 46N10

\section{Introduction}

Given a locally convex Hausdorff topological vector space $X$ and an
infinite family $\left( f_{i}\right) _{i\in I}\subset \left( \mathbb{R}%
_{\infty }\right) ^{X},$ where $\mathbb{R}_{\infty }:=\mathbb{R\cup }\left\{
+\infty \right\} ,$ of objective proper functions, we are concerned with the
uncertain problem of minimizing a finite but unknown sum of the objective
functions $f_{i}.$ Adopting the robust optimization approach under
uncertainty (see \cite{BJL13}, \cite{DGLV18}, \cite{DGLV18B}, \cite{LJL11}),
and taking the set $\mathcal{F}\left( I\right) $ of non-empty finite subsets
of $I$ as uncertainty set, the robust counterpart of this uncertain problem
is%
\begin{equation}
\mathrm{\left( RP\right)} \ \ \ \ \ \ \ \inf\limits_{x\in
X}\sup\limits_{J\in \mathcal{F}\left( I\right) }\sum\nolimits_{i\in J}{{f_{i}%
}}\left( x\right) .\text{ }  \label{1.1}
\end{equation}

This kind of problem arises in situations where one must minimize the
\textit{robust-}$L_{p}$ \textit{(pseudo)\ norm} function $\Vert h \Vert_p$
defined on $X$ by
\begin{equation*}
\Vert h\Vert _{p} (x):=\left[ \sup\limits_{J\in \mathcal{F}\left( I\right)
}\sum\limits_{i\in J}\left\vert h{{_{i}}}(x) \right\vert ^{p}\right] ^{\frac{%
1}{p}} = \sup\limits_{J\in \mathcal{F}\left( I\right) } \left[
\sum\limits_{i\in J}\left\vert h{{_{i}}}(x) \right\vert ^p \right]^{\frac{1}{%
p}} \in \mathbb{R}_{ \infty},
\end{equation*}%
for any $h = (h_i)_{i \in I} \subset \mathbb{R}^X$ and $p \geq 1$. Since the
exact value of $\Vert h\Vert_{p}$ can hardly be computed, in practice it
should be replaced by the maximum of $\left[ \sum\nolimits_{i\in
J}\left\vert h{{_{i}}}\right\vert ^{p}\right] ^{\frac{1}{p}}$ for a sample
of non-empty finite sets $J$ picked at random from $I.$ This way, $\left[
\sum\nolimits_{i\in J}\left\vert h{{_{i}}} \right\vert ^{p}\right] ^{\frac{1%
}{p}}$ can be interpreted as an uncertain function with uncertain parameter $%
J$ ranging on the uncertainty set $\mathcal{F}\left( I\right) .$

As a first example, in the extension of the least squares linear regression
model to the case of infinite point clouds $\left\{ (t_{i},s_{i}),\,i\in
I\right\} \subset \mathbb{R}^{2}$, when the shape of the latter set suggests
a linear dependence of magnitude $s$ with respect to magnitude $t,$ the
problem consists in computing the ordinate at the origin, $x_{1},$ and the
slope, $x_{2},$ of the line $s=x_{1}+x_{2}t$ better fitted to that set. To
find it, one should minimize $\Vert h\Vert_{2}^{2} \left( x_{1},x_{2}\right)
$ \ on $\mathbb{R}^{2},$\ where the $i-$th component of the residual
function $h,$ is ${h{_{i}}}\left( x_{1},x_{2}\right)
:=x_{1}+x_{2}t_{i}-s_{i},$ $i\in I.$ In the terminology of robust
optimization, the uncertain objective function $\sum\nolimits_{i\in
J}\left\vert h{{_{i}}}\left( x\right) \right\vert ^{2}$ at $x$\ represents
the sum of squares error for the line $s=x_{1}+x_{2}t$ relative to the
finite point cloud $\left\{ (t_{i},s_{i}),\,i\in J\right\} ,$ the worst-case
objective function $\sup\nolimits_{J\in \mathcal{F}\left( I\right)
}\sum\nolimits_{i\in J}\left\vert h{{_{i}}}\left( x\right) \right\vert ^{2}$
is the least upper bound, for $J\in \mathcal{F}\left( I\right) ,$ of the
errors corresponding to that line, and a robust optimal solution of $\left(%
\mathrm{RP}\right) $ is a best infinite regression line for the point cloud $%
\left\{ (t_{i},s_{i}),\,i\in I\right\} $.

A second example comes from the search of a best approximate solution to an
inconsistent system $\left\{ \left\langle a_{i},x\right\rangle \leq b_{i},\
i\in I\right\} $ in $\mathbb{R}^{n}$. Denote by $h\left( x\right) $ the
residual of$\ x\in \mathbb{R}^{n},$ i.e., $h_{i}\left( x\right) :=\max
\left\{ \left\langle a_{i},x\right\rangle -b_{i}.0\right\} ,$ $i\in I.$
Assuming that $I$ is the union of a discrete set with a finite union of
pairwise disjoint boxes as well as the continuity on these boxes of the
function $i\longmapsto \left( a_{i},b_{i}\right) ,$ \cite{GHL18} analyzes
the minimization of the components of the residual function $h$, involving
integrals whose existence is guaranteed by the continuity assumption. One
can get rid of any assumption on $I$ and the function $i\longmapsto \left(
a_{i},b_{i}\right) $ by considering the minimization of the robust
pseudonorm function $\Vert h\left( x\right) \Vert _{p}$ for an arbitrary
infinite set $I,$ in which case an optimal solution of $\mathrm{\left(
RP\right) }$ provides a best robust-$L_{p}$ approximate solution of $\left\{
\left\langle a_{i},x\right\rangle \leq b_{i},\ t\in I\right\} .$

The third example involving robust sum functions, not related with the above
robust-$L_{p}$ norm, is inspired in the classic portfolio model for a finite
set $I$ of assets with expected return $r_{i}$ and estimated covariance $%
v_{ij}$ of the returns of assets $i,j\in I:$%
\begin{equation*}
\begin{array}{lll}
\mathrm{\left( P\right) } & \max  & \sum\limits_{i\in I}r_{i}x_{i} \\
& \min  & \sum\limits_{i,j\in I}v_{ij}x_{i}x_{j} \\
& \text{s.t.} & \sum\limits_{i\in I}x_{i}=1, \\
&  & x_{i}\geq 0,i\in I,%
\end{array}%
\end{equation*}%
where the two objectives consist in maximizing the expected return of the
portfolio and minimizing its volatility (identified here with the risk of
the portfolio). Taking into account the almost unlimited number of existing
assets in the global economy, it is natural to replace $I$ by $\mathbb{N}$
in $\left( \mathrm{P}\right) ,$ the decision space $\mathbb{R}^{\left\vert
I\right\vert }$ by $X:=\mathbb{R}^{\mathbb{N}},$ the first objective by the
minimization of $f\left( x\right) :=\sup\nolimits_{J\in \mathcal{F}\left(
\mathbb{N}\right) }\sum\nolimits_{i\in J}{{f_{i}}}\left( x\right) ,$ with ${{%
f_{i}}}\left( x\right) :=-r_{i}x_{i},$ the second objective function by $%
g\left( x\right) :=\sup\nolimits_{K\in \mathcal{F}\left( \mathbb{N}%
^{2}\right) }\sum\nolimits_{\left( i,j\right) \in K}{g{_{ij}}}\left(
x\right) ,$ with ${g{_{ij}}}\left( x\right) =\sum\nolimits_{\left(
i,j\right) \in K}v_{ij}x_{i}x_{j},$ and the first constraint by $h\left(
x\right) =1,$ where $h$ is the robust sum $h\left( x\right)
:=\sup\nolimits_{J\in \mathcal{F}\left( \mathbb{N}\right)
}\sum\nolimits_{i\in J}x_{i},$ giving rise to a bi-objective infinite
dimensional optimization problem $\mathrm{\left( P_{1}\right) }$ involving
robust sums of linear functions and quadratic forms. When the decision maker
is able to fix a volatility threshold $\sigma ^{2},$ the problem to be
solved is a scalar one involving robust sum functions:%
\begin{equation*}
\begin{array}{lll}
\mathrm{\left( P_{2}\right) } & \max  & f\left( x\right)  \\
& \text{s.t.} & g\left( x\right) \leq \sigma ^{2}, \\
&  & h\left( x\right) =1, \\
&  & x_{i}\geq 0,i\in I.%
\end{array}%
\end{equation*}

The aim of this paper is to establish some duality principles for the
problem $\mathrm{\left( RP\right) }$ and to characterize in various ways the
zero duality gap property. We call \textit{robust sum} of the family $\left(
f_{i}\right) _{i\in I} \subset (\mathbb{R}_\infty)^X$, represented by $%
\sum\nolimits_{i\in I}^{R} f_{i}: X\longrightarrow \mathbb{R}_\infty$, the
objective function of $\mathrm{\left( RP\right) }$, namely,%
\begin{equation*}
\sum\nolimits_{i\in I}^{R}{{f_{i}}\left( x\right) :=}\sup\limits_{J\in
\mathcal{F}\left( I\right) }\sum\limits_{i\in J}{{f_{i}}}\left( x\right)
,\forall x\in X.
\end{equation*}
The term "robust sum" is not new in the literature, but it has been only
used in the framework of the uncertain optimization of finite sums (see,
e.g., \cite{AD}, \cite{DN}).

In the case where all functions $f_{i}$ are non-negative,
\begin{equation}
\sum\nolimits_{i\in I}^{R}{{f_{i}}\left( x\right) =}\sum\limits_{i\in I}{{%
f_{i}}}\left( x\right) :=\lim\limits_{J\in \mathcal{F}\left( I\right)
}\sum\limits_{i\in J}{{f_{i}}}\left( x\right) ,\forall x\in X,  \label{1.3}
\end{equation}%
where the limit is taken respect to the directed set $\mathcal{F}\left(
I\right) $ ordered by the inclusion relation. The advantage of the robust
sum $\sum\nolimits_{i\in I}^{R}{{f_{i}}}$ in comparison with the \textit{%
infinite sum} $\sum\nolimits_{i\in I}{{f_{i}}}$ is that $\sum\nolimits_{i\in
I}^{R}{{f_{i}}}\left( x\right) $ is well defined for each $x\in X$ while $%
\sum\nolimits_{i\in I}{{f_{i}}}\left( x\right) $ may not exist (see Remark %
\ref{Rem2.1} and Lemma \ref{Lemma 2.6} below). Formulas for the
subdifferential of $\sum\nolimits_{i\in I}{{f_{i}}}$ in the case that all
functions $f_{i}$ are continuous have been given in \cite{Zheng98} and \cite[%
Proposition 2.3]{ZN04}, while duality theorems on infinite sums of proper,
convex and lower semicontinuous (lsc in short) functions can be found in
\cite[Section 3]{LN08}. The mentioned subdifferential formulas and duality
theorems for $\sum\nolimits_{i\in I}{{f_{i}}}$ have been used in \cite[%
Proposition 2.3]{ZN04} and \cite[Section 5]{LN08} to obtain error bounds for
convex infinite systems and optimality conditions for convex infinite
programs, respectively.

Throughout this paper we assume that all functions $f_{i},$ $i\in I,$ are
proper, as well as their robust sum $f:=\sum\nolimits_{i\in I}^{R}{{f_{i}.}}$
The paper is organized as follows. Section 2 introduces the robust sum of an
infinite family in $\mathbb{R}_{\infty }$ and analyzes its relationship with
the infinite sum of the family. Sections 3, 4 and 5 provide results
characterizing weak duality, zero duality gap, and strong duality, for the
robust sum of a family of arbitrary functions, respectively, in terms of
multifunctions associated with $\left( f_{i}\right) _{i\in I}$. Section 6
analyzes the robust sum under the assumption that $\left( f_{i}\right)
_{i\in I}$ is a family of proper, lsc and convex functions; the main result
of this section is Theorem \ref{Th5}, which characterizes the strong zero
duality gap of $f$ under a closedness assumption instead of $\varepsilon $%
-subdifferentials and epigraphs of the family of corresponding conjugate
functions, as in \cite[Theorem 3.2]{LN08} for $\sum\nolimits_{i\in I}{{f_{i}}%
}.$ Finally, Section 7 provides a stable zero duality theorem for the
infinite sum of proper, lsc, and non-negative convex functions (as in the
above infinite regression problem).

\section{Some rules for the robust sum}

We associate with a given infinite family $\left( a_{i}\right) _{i\in
I}\subset \mathbb{R}_{\infty }$ its \textit{robust sum}
\begin{equation}
\sum\nolimits_{i\in I}^{R}a_{i}:=\sup\limits_{J\in \mathcal{F}\left(
I\right) }\sum\limits_{i\in J}a_{i},  \label{0.1}
\end{equation}%
together with its \textit{inferior and superior limits},
\begin{equation*}
\liminf\limits_{J\in \mathcal{F}\left( I\right) }\sum\limits_{i\in
J}a_{i}:=\sup\limits_{J\in \mathcal{F}\left( I\right) }\inf_{J\subset K\in
\mathcal{F}\left( I\right) }\sum\limits_{i\in K}a_{i},
\end{equation*}%
and%
\begin{equation*}
\limsup\limits_{J\in \mathcal{F}\left( I\right) }\sum\limits_{i\in
J}a_{i}:=\inf\limits_{J\in \mathcal{F}\left( I\right) }\sup\limits_{J\subset
K\in \mathcal{F}\left( I\right) }\sum\limits_{i\in K}a_{i},
\end{equation*}%
respectively.

\begin{lemma}
\label{Lemma 2.1} One has%
\begin{equation}
-\infty <\sup\limits_{i\in I}a_{i}\leq \sum\nolimits_{i\in I}^{R}a_{i}\leq
+\infty ,  \label{0.2}
\end{equation}%
and
\begin{equation}
-\infty <\liminf\limits_{J\in \mathcal{F}\left( I\right) }\sum\limits_{i\in
J}a_{i}\leq \sum\nolimits_{i\in I}^{R}a_{i}\leq +\infty .  \label{0.3}
\end{equation}
\end{lemma}

\textit{Proof.} Let $j\in I.$ Setting $J=\left\{ j\right\} $ in (\ref{0.1})
we get%
\begin{equation*}
-\infty <a_{j}\leq \sum\nolimits_{i\in I}^{R}a_{i}\leq +\infty .
\end{equation*}%
Taking the supremum over $j\in I,$ (\ref{0.2}) holds true.

Let $J\in \mathcal{F}\left( I\right) .$ We have%
\begin{equation*}
-\infty \leq \inf_{ K\in \mathcal{F}\left( I\right) }\sum\limits_{i\in
K}a_{i}\leq \sum\limits_{i\in J}a_{i}\leq \sum\nolimits_{i\in
I}^{R}a_{i}\leq +\infty .
\end{equation*}%
Taking the supremum over $J\in \mathcal{F}\left( I\right) ,$ (\ref{0.3})
holds true. \qquad \qed

We also define the \textit{infinite sum} of the family $\left( a_{i}\right)
_{i\in I}$ as
\begin{equation*}
\sum\limits_{i\in I}a_{i}:=\lim\limits_{J\in \mathcal{F}\left( I\right)
}\sum\limits_{i\in J}a_{i},
\end{equation*}%
provided that the unconditional limit $\lim\limits_{J\in \mathcal{F}\left(
I\right) }\sum\nolimits_{i\in J}a_{i}$ exists as a member of $\overline{%
\mathbb{R}},$ i.e.,%
\begin{equation*}
-\infty \leq \liminf\limits_{J\in \mathcal{F}\left( I\right)
}\sum\limits_{i\in J}a_{i}=\limsup\limits_{J\in \mathcal{F}\left( I\right)
}\sum\limits_{i\in J}a_{i}\leq +\infty .
\end{equation*}

In the case when $\left( a_{i}\right) _{i\in I}\in \left[ 0,+\infty \right]
, $ we have
\begin{equation*}
0\leq \sum\nolimits_{i\in I}^{R}a_{i}=\sum\limits_{i\in I}a_{i}\leq +\infty .
\end{equation*}

For each $\theta \in \overline{\mathbb{R}}$ we consider $\theta ^{+}:=\max
\left\{ \theta ,0\right\} $ and $\theta ^{-}:=\max \left\{ -\theta
,0\right\} .$

\begin{lemma}
\label{Lemma 2.2} $\left( \sum\nolimits_{i\in I}^{R}a_{i}\right)
^{+}=\sum\nolimits_{i\in I}^{R}a_{i}^{+}=\sum\limits_{i\in I}a_{i}^{+}.$
\end{lemma}

\textit{Proof.} Since $a_{i}\leq a_{i}^{+}$ and $a_{i}^{+}\geq 0$ for all $%
i\in I,$ we have%
\begin{equation*}
\sum\nolimits_{i\in I}^{R}a_{i}\leq \sum\nolimits_{i\in
I}^{R}a_{i}^{+}=\sum\limits_{i\in I}a_{i}^{+}.
\end{equation*}%
Since $\sum\nolimits_{i\in I}a_{i}^{+}\geq 0$ we obtain $\left(
\sum\nolimits_{i\in I}^{R}a_{i}\right) ^{+}\leq \sum\nolimits_{i\in
I}a_{i}^{+}.$ Let us prove the reverse inequality. Let $J\in \mathcal{F}%
\left( I\right)$ and $K_J:=\left\{ i\in J:a_{i}>0\right\} .$ If $%
K_J=\emptyset $ then $\sum\nolimits_{i\in J}a_{i}^{+}=0 \leq \left(
\sum\nolimits_{i\in I}^{R}a_{i}\right) ^{+}$. If, alternatively, $K_{J} \neq
\emptyset $ then
\begin{equation*}
0\leq \sum\limits_{i\in J}a_{i}^{+}=\sum\limits_{i\in K_{J}}a_{i}\leq
\sum\nolimits_{i\in I}^{R}a_{i}\leq \left( \sum\nolimits_{i\in
I}^{R}a_{i}\right) ^{+}.
\end{equation*}%
In both cases we have $\sum\nolimits_{i\in J}a_{i}^{+}\leq \left(
\sum\nolimits_{i\in I}^{R}a_{i}\right) ^{+},$ and, since $J\in \mathcal{F}%
\left( I\right) $ is arbitrary, we obtain%
\begin{equation*}
\sum\limits_{i\in I}a_{i}^{+}=\sum\nolimits_{i\in I}^{R}a_{i}^{+}\leq \left(
\sum\nolimits_{i\in I}^{R}a_{i}\right) ^{+},
\end{equation*}%
and the proof is complete. \qquad \qed

As an immediate consequence of Lemma \ref{Lemma 2.2} we have:

\begin{lemma}
\label{Lemma 2.3} One has
\begin{equation*}
\sum\nolimits_{i\in I}^{R}a_{i}\in \mathbb{R\Longleftrightarrow }%
\sum\limits_{i\in I}a_{i}^{+}<+\infty .
\end{equation*}
\end{lemma}

\begin{lemma}
\label{Lemma 2.4} Next statements are equivalent:\newline
$(i)$ $\sup\limits_{i\in I}a_{i}\geq 0.$\newline
$(ii)$ $\sum\nolimits_{i\in I}^{R}a_{i}\geq 0.$\newline
$(iii)$ $\sum\nolimits_{i\in I}^{R}a_{i}=\sum\nolimits_{i\in I}a_{i}^{+}.$
\end{lemma}

\textit{Proof.} One has $\left[ (i)\Longrightarrow (ii)\right] $ by Lemma %
\ref{Lemma 2.1} and $\left[ (ii)\Longrightarrow (i)\right] $ by Lemma \ref%
{Lemma 2.2}. Assume now that $(i)$ does not hold, i.e., there exists $%
\varepsilon >0$ such that $a_{i}\leq -\varepsilon $ for all $i\in I.$ For
each $J\in \mathcal{F}\left( I\right) $ we have $\sum\nolimits_{i\in
J}a_{i}\leq -\varepsilon \times \func{card}J\leq -\varepsilon .$ Since $J$
is arbitrary we get $\sum\nolimits_{i\in I}^{R}a_{i}\leq -\varepsilon $ and $%
(iii)$ does not hold. So, $\left[ (iii)\Longrightarrow (i)\right] $ and the
proof is complete. \qquad \qed

\begin{lemma}
\label{Lemma 2.5} One has%
\begin{equation}
\sum\nolimits_{i\in I}^{R}a_{i}=\left\{
\begin{array}{ll}
\sum\nolimits_{i\in I}a_{i}^{+}, & \text{if } \ \ \sup\limits_{i\in
I}a_{i}\geq 0, \\
\sup\limits_{i\in I}a_{i}, & \text{if } \ \ \sup\limits_{i\in I}a_{i}\leq 0.%
\end{array}%
\right.  \label{0.4}
\end{equation}
\end{lemma}

\textit{Proof.} If $\sup\limits_{i\in I}a_{i}\geq 0$, (\ref{0.4}) holds by
Lemma \ref{Lemma 2.4}. Assume now that $\sup\limits_{i\in I}a_{i}\leq 0.$ By
Lemma \ref{Lemma 2.1} we have just to check that $\sum\nolimits_{i\in
I}^{R}a_{i}\leq \sup\limits_{i\in I}a_{i}.$ Let $J\in \mathcal{F}\left(
I\right) .$ Picking $j\in J$, we have
\begin{equation*}
\sum\limits_{i\in J}a_{i}\leq a_{j}\leq \sup\limits_{i\in I}a_{i},
\end{equation*}%
and, since $J$ is arbitrary, we are done. \qquad \qed

\begin{remark}
\label{Rem2.1} We note that $\sum\nolimits_{i\in I}^{R}a_{i}$ always exists
in $\mathbb{R}_{\infty }$ while $\sum\nolimits_{i\in I}a_{i}$ may not exist
in $\overline{\mathbb{R}}.$ This is for instance the case when $I=\mathbb{N}$
and $a_{i}=\left( -1\right) ^{i}$ or $a_{i}=\frac{\left( -1\right) ^{i}}{i}.$
In both cases we have $\sum\nolimits_{i\in I}a_{i}^{+}=+\infty $ and, by
Lemma \ref{Lemma 2.3}, $\sum\nolimits_{i\in I}^{R}a_{i}=+\infty ,$ while
there are subnets of $\left\{ \sum\nolimits_{i\in J}a_{i}\right\} _{J\in
\mathcal{F}\left( I\right) }$ converging towards distinct limits, so that $%
\sum\nolimits_{i\in I}a_{i}$ does not exist in $\overline{\mathbb{R}}.$
However, in the case when $\sum\nolimits_{i\in I}^{R}a_{i}\in \mathbb{R},$ $%
\sum\nolimits_{i\in I}a_{i}$ does exist in $\mathbb{R\cup }\left\{ -\infty
\right\} $ as the next lemma shows.
\end{remark}

\begin{lemma}
\label{Lemma 2.6} Assume that $\sum\nolimits_{i\in I}^{R}a_{i}\in \mathbb{R}%
. $ Then either $\sum\nolimits_{i\in I}a_{i}^{-}<+\infty $ or $%
\sum\nolimits_{i\in I}a_{i}^{-}=+\infty$, and in the first case it holds $%
\sum\nolimits_{i\in I}a_{i}\in \mathbb{R}$ while in the latter one, $%
\sum\nolimits_{i\in I}a_{i}=-\infty .$
\end{lemma}

\textit{Proof.} Since $\sum\nolimits_{i\in I}^{R}a_{i}\in \mathbb{R},$ Lemma %
\ref{Lemma 2.3} says that $\sum\nolimits_{i\in I}a_{i}^{+}<+\infty .$

Assume that $\sum\nolimits_{i\in I}a_{i}^{-}<+\infty .$ Then,
\begin{equation*}
\sum\nolimits_{i\in I}a_{i}^{+}-\sum\nolimits_{i\in
I}a_{i}^{-}=\sum\nolimits_{i\in I}\left( a_{i}^{+}-a_{i}^{-}\right)
=\sum\nolimits_{i\in I}a_{i}\in \mathbb{R}.
\end{equation*}

Assume that $\sum\nolimits_{i\in I}a_{i}^{-}=+\infty $ and let us prove that
$\sum\nolimits_{i\in I}a_{i}=-\infty .$

Let $r\in \mathbb{R}$ and $s:=\sum\nolimits_{i\in I}a_{i}^{+}\in \mathbb{R}$%
. There exist $J_{1},J_{2}\in \mathcal{F}\left( I\right) $ such that,
\begin{equation*}
\forall J\in \mathcal{F}\left( I\right), \ J_{1}\subset J \ \
\Longrightarrow \ \ \sum\nolimits_{i\in J}a_{i}^{+}\leq s+1 ,
\end{equation*}%
\begin{equation*}
\forall J\in \mathcal{F}\left( I\right), \ J_{2}\subset J\Longrightarrow
\sum\nolimits_{i\in J}a_{i}^{-}\geq s+1-r.
\end{equation*}
Thus, for each $J\in \mathcal{F}\left( I\right) $ such that $J_{1}\cup
J_{2}\subset J$, we have%
\begin{equation*}
\sum\nolimits_{i\in J}a_{i}=\sum\nolimits_{i\in
J}a_{i}^{+}-\sum\nolimits_{i\in J}a_{i}^{-}\leq r,
\end{equation*}%
which means that $\sum\nolimits_{i\in I}a_{i}=-\infty .$ \qquad \qed

\begin{example}
\label{Exam1} Let $I=\mathbb{N}$ and, for each $i\in I,$
\begin{equation*}
a_{i}=\left\{
\begin{array}{ll}
\frac{1}{i^{2}}, & \text{if }i\text{ is even}, \\
-\frac{1}{i}, & \text{if }i\text{ is odd}.%
\end{array}%
\right.
\end{equation*}%
By Lemma \ref{Lemma 2.5} we have $\sum\nolimits_{i\in
I}^{R}a_{i}=\sum\nolimits_{i\in I}a_{i}^{+}=\frac{\pi ^{2}}{24}$ and, since $%
\sum\nolimits_{i\in I}a_{i}^{-}=+\infty ,$ by Lemma \ref{Lemma 2.6}, we have
$\sum\nolimits_{i\in I}a_{i}=-\infty .$
\end{example}

\section{Weak duality}

We now introduce the notation that will be used in the rest of the paper.
The topological dual space of $X$ is denoted by $X^{\ast }.$ We denote by $%
0_{X}$\ and $0_{X}^{\ast }$ the null vector of $X$ and $X^{\ast },$
respectively. The closure of a subset $A \subset X$ will be denoted by $%
\overline{A}$ and the same symbol will be used for the closure of a subset
of the dual space $X^\ast$.

Given a function $h\in \overline{\mathbb{R}}^{X},$ its domain is the set $%
\limfunc{dom}h:=\{x\in X:h(x)<+\infty \},$ its epigraph is $\limfunc{epi}%
h:=\{\left( x,r\right) \in X\times \mathbb{R}:h(x)\leq r\},$ its strict
epigraph is $\limfunc{epi}\nolimits_{s}h:=\{\left( x,r\right) \in X\times
\mathbb{R}:h(x)<r\},$ and its Fenchel conjugate is the function $h^{\ast
}\in \overline{\mathbb{R}}^{X^{\ast }}$ such that $h^{\ast }(x^{\ast
}):=\sup \{\langle x^{\ast },x\rangle -h(x):x\in X\}$ for any $x^{\ast }\in
X^{\ast }$. Moreover, the lsc hull of $h$ is the function $\overline{h}\in
\overline{\mathbb{R}}^{X}$ whose epigraph $\limfunc{epi}\overline{h}$ is the
closure of $\limfunc{epi}h$ in $X\times \mathbb{R}.$

Given $\varepsilon \in \mathbb{R}$, we denote by $\left[ h\leq \varepsilon %
\right] :=\left\{ x\in X:h(x)\leq \varepsilon \right\} $\ the lower level
set of $h$ at level $\varepsilon .$ The definition of the strict lower level
set $\left[ h<\varepsilon \right] $ is similar.

Given $\varepsilon \geq 0,$ we define the $\varepsilon -$minimizers of $h$
as
\begin{equation*}
\varepsilon -\func{argmin}h:=\left\{
\begin{array}{ll}
\left\{ x\in X:h\left( x\right) \leq \inf\limits_{X}h+\varepsilon \right\} ,
& \text{if }\inf\limits_{X}h\in \mathbb{R}, \\
\emptyset , & \text{else.}%
\end{array}%
\right.
\end{equation*}

Given $a\in X$ and $\varepsilon \geq 0$, we denote by%
\begin{equation*}
\partial ^{\varepsilon }h(a):=\left\{
\begin{array}{ll}
\{x^{\ast }\in X^{\ast }\,:\,h(x)\geq h(a)+\langle x^{\ast },x-a\rangle
-\varepsilon ,\forall x\in X\}, & \text{if }h(a)\in \mathbb{R}, \\
\emptyset , & \text{else,}%
\end{array}%
\right.
\end{equation*}%
the $\varepsilon -$subdifferential of $h$ at $a.$ For $\varepsilon =0$ one
sets $\partial h(a)$ instead of $\partial ^{0}h(a).$ By definition, $%
\partial ^{\varepsilon }h:X\rightrightarrows X^{\ast }$ is a multifunction
whose inverse multifunction we denote by $M^{\varepsilon }h:X^{\ast
}\rightrightarrows X.$ For each $x^{\ast }\in X^{\ast }$ one has
\begin{equation*}
M^{\varepsilon }h(x^{\ast })=\left\{
\begin{array}{ll}
\!\varepsilon -\limfunc{argmin}(h-x^{\ast }), & \!\!\!\mathrm{if}\ h^{\ast
}(x^{\ast })\in \mathbb{R}, \\
\emptyset , & \!\!\!\text{\textrm{else}}.%
\end{array}%
\right.
\end{equation*}%
The multifunction $M^{\varepsilon }h(\cdot)$ will be of a crucial importance
in the paper. Notice that, with the rule $\left( +\infty \right) -\left(
-\infty \right) =\left( -\infty \right) +\left( +\infty \right) =+\infty ,$
one has%
\begin{equation}  \label{eqsubdif}
x\in M^{\varepsilon }h(x^{\ast })\Longleftrightarrow x^{\ast }\in \partial
^{\varepsilon }h(x)\Longleftrightarrow h(x)+h^{\ast }(x^{\ast })\leq \langle
x^{\ast },x\rangle +\varepsilon.
\end{equation}

We are now turning back to the problem $\left(\mathrm{RP}\right) $ defined
in \eqref{1.1} by an infinite family $\left( f_{i}\right) _{i\in I}\subset
\left( \mathbb{R}_{\infty }\right) ^{X}$ of proper functions with $f
=\sum\nolimits_{i\in I}^{R}{{f_{i}}}$, which is assumed to be proper as
well. Note that as the functions $f$, $f_{i}$ are proper, the conjugate
functions $f^\ast$, $f_{i}^{\ast },$ $i\in I,$ never take the value $-\infty
.$

For each $x^*\in X^*$ consider the dual pair of problems
\begin{align*}
(\mathrm{RP}_{x^*})\quad & \inf _{x\in X} [f(x)-\langle x^*, x\rangle], \\
(\mathrm{RD}_{x^*})\quad & \sup_{\substack{ I\in \mathcal{F}(I)  \\ %
(x^*_i)_{i\in I} \in (X^*)^I  \\ \sum_{i\in J} x^*_i=x^*}} -\sum_{i\in J}
f^*_i(x^*_i).
\end{align*}
It is clear that (RP) is nothing else but $(\mathrm{RP}_{0_{X^*}})$ and from
now on, we will write (RP) and (RD) instead of $(\mathrm{RP}_{0_{X^*}})$ and
$(\mathrm{RD}_{0_{X^*}})$, respectively. Note that (RD) is nothing but the
optimistic dual problem of (RP).

Let us now introduce the function $\varphi :X^{\ast }\longrightarrow
\overline{\mathbb{R}}$ defined as%
\begin{equation*}
\varphi \left( x^{\ast }\right) :=\inf_{J\in \mathcal{F}\left( I\right)
}\left\{ \dsum\limits_{i\in J}{{f_{i}^{\ast }}}\left( x_{i}^{\ast }\right)
:\left( x_{i}^{\ast }\right) _{i\in J}\in \left( X^{\ast }\right)
^{J},\dsum\limits_{i\in J}x_{i}^{\ast }=x^{\ast }\right\} ,\forall x^{\ast
}\in X^{\ast }.
\end{equation*}
Then it is clear that for each $x^* \in X^\ast$,
\begin{equation*}
\inf (\mathrm{RP}_{x^*}) = - f^\ast (x^*) \ \ \mathrm{and} \ \ \sup (\mathrm{%
RD}_{x^*}) = - \varphi (x^*) .
\end{equation*}

\begin{proposition}[Weak duality]
\label{Prop1}For each $x^{\ast }\in X^{\ast }$ we have%
\begin{equation*}
-\infty \leq \sup (\mathrm{RD}_{x^*}) \leq \inf (\mathrm{RP}_{x^*}) <
+\infty ,
\end{equation*}
or, equivalently,
\begin{equation*}
-\infty <{{f^{\ast }}}\left( x^{\ast }\right) \leq \varphi \left( x^{\ast
}\right) \leq +\infty .
\end{equation*}
\end{proposition}

\textit{Proof.} Since $f$ is proper, its conjugate does not take the value $%
-\infty .$ Let $x^{\ast }\in X^{\ast },$ $J\in \mathcal{F}\left( I\right) ,$
$\left( x_{i}^{\ast }\right) _{i\in J}\in \left( X^{\ast }\right)
^{J},\dsum\limits_{i\in J}x_{i}^{\ast }=x^{\ast },$ and $x\in X.$ One has to
check that $\left\langle x^{\ast },x\right\rangle -f\left( x\right) \leq
\dsum\limits_{i\in J}{{f_{i}^{\ast }}}\left( x_{i}^{\ast }\right) .$ By
definition of ${{f_{i}^{\ast }}}$ we have%
\begin{equation*}
\dsum\limits_{i\in J}{{f_{i}^{\ast }}}\left( x_{i}^{\ast }\right) \geq
\dsum\limits_{i\in J}\left( \left\langle x_{i}^{\ast },x\right\rangle -
f_i\left( x\right) \right) =\left\langle x^{\ast },x\right\rangle
-\dsum\limits_{i\in J}{{f_{i}}}\left( x\right) \geq \left\langle x^{\ast
},x\right\rangle - f(x),
\end{equation*}%
as by the definition of $f$, $\dsum\nolimits_{i\in J}{{f_{i}}}\left(
x\right) \leq f\left( x\right) $ for all $x\in X$. \qquad \qed

\section{Zero duality gap}

\begin{definition}
We say that the robust sum problem $(\mathrm{RP}_{x^*}) $ has zero duality
gap at a given $x^{\ast }\in X^{\ast }$ if
\begin{equation}
\inf (\mathrm{RP}_{x^*}) = \sup (\mathrm{RD}_{x^*}),  \label{2.1}
\end{equation}
(or equivalently, ${{f^{\ast }}}\left( x^{\ast }\right) =\varphi \left(
x^{\ast }\right)) $. If \ (\ref{2.1}) holds at each $x^{\ast }\in X^{\ast },$
we will say that the problem $(\mathrm{RP}_{x^*}) $ has stable zero duality
gap.
\end{definition}

The characterization of the zero duality gap involves two mutually inverse
multifunctions associated with the given family of functions $\left(
f_{i}\right) _{i\in I}$ with robust sum $f.$

For each $\alpha \geq 0$ let us define $S_{f}^{\alpha }:X\rightrightarrows
\mathcal{F}\left( I\right) $ such that%
\begin{equation*}
S_{f}^{\alpha }\left( x\right) :=\left\{
\begin{array}{ll}
\left\{ J\in \mathcal{F}\left( I\right) :f\left( x\right) \leq
\sum\nolimits_{i\in J}{{f_{i}}}\left( x\right) +\alpha \right\} , & \text{if
}x\in \limfunc{dom}f, \\
\emptyset , & \text{else,}%
\end{array}%
\right.
\end{equation*}%
and $T_{f}^{\alpha }:\mathcal{F}\left( I\right) \rightrightarrows X$ such
that%
\begin{equation*}
T_{f}^{\alpha }\left( J\right) :=\left\{ x\in \limfunc{dom}f:f\left(
x\right) \leq \sum\nolimits_{i\in J}{{f_{i}}}\left( x\right) +\alpha
\right\} .
\end{equation*}%
So, for any $\left( x,J\right) \in X\times \mathcal{F}\left( I\right) $ we
have%
\begin{equation*}
J\in S_{f}^{\alpha }\left( x\right) \Longleftrightarrow x\in T_{f}^{\alpha
}\left( J\right) .
\end{equation*}

Let us now put in light a necessary condition for the robust sum problem to
have zero duality gap at a given $x^{\ast }\in X^{\ast }.$ So, assume that%
\begin{equation}
{{f^{\ast }}}\left( x^{\ast }\right) \geq \varphi \left( x^{\ast }\right)
\in \mathbb{R}  \label{1}
\end{equation}%
and let $\varepsilon \geq 0.$ For any $\eta >0$ and $x\in M^{\varepsilon
}f(x^{\ast }),$ one has%
\begin{equation*}
f\left( x\right) -\left\langle x^{\ast },x\right\rangle \leq -{{f^{\ast }}}%
\left( x^{\ast }\right) +\varepsilon <-\varphi \left( x^{\ast }\right)
+\varepsilon +\eta .
\end{equation*}%
Then, by definition of $\varphi ,$ there exist $J\in \mathcal{F}\left(
I\right) $ and $\left( x_{i}^{\ast }\right) _{i\in J}\in \left( X^{\ast
}\right) ^{J}$ such that $\dsum\nolimits_{i\in J}x_{i}^{\ast }=x^{\ast }$ and%
\begin{equation*}
f\left( x\right) -\left\langle x^{\ast },x\right\rangle \leq
-\dsum\limits_{i\in J}{{f_{i}^{\ast }}}\left( x_{i}^{\ast }\right)
+\varepsilon +\eta .
\end{equation*}%
This last inequality can be rewritten as%
\begin{equation*}
\left[ f\left( x\right) -\sum\limits_{i\in J}{{f_{i}}}\left( x\right) \right]
+\dsum\limits_{i\in J}\left[ {{f_{i}}}\left( x\right) +{{f_{i}^{\ast }}}%
\left( x_{i}^{\ast }\right) -\left\langle x_{i}^{\ast },x\right\rangle %
\right] \leq \varepsilon +\eta .
\end{equation*}%
All the above brackets being non-negative, there exists $\left( \alpha
,\left( \varepsilon _{i}\right) _{i\in J}\right) \in \mathbb{R}_{+}\times
\mathbb{R}_{+}^{J}$ such that $\alpha +\sum\nolimits_{i\in J}\varepsilon
_{i}=\varepsilon +\eta $ and
\begin{equation*}
f\left( x\right) -\sum\limits_{i\in J}{{f_{i}}}\left( x\right) \leq \alpha
\text{ and }{{f_{i}}}\left( x\right) +{{f_{i}^{\ast }}}\left( x_{i}^{\ast
}\right) -\left\langle x_{i}^{\ast },x\right\rangle \leq \varepsilon
_{i},i\in J.
\end{equation*}%
In other words,
\begin{equation}
x\in T_{f}^{\alpha }\left( J\right) \text{ and }x\in M^{\varepsilon
_{i}}f(x_{i}^{\ast }),\forall i\in J.  \label{2}
\end{equation}

Hence we have quoted that for any $x\in M^{\varepsilon }f(x^{\ast })$ and
any $\eta >0,$ there exist $J\in \mathcal{F}\left( I\right) ,$ $\left(
x_{i}^{\ast }\right) _{i\in J}\in \left( X^{\ast }\right) ^{J}$ and $\left(
\alpha ,\left( \varepsilon _{i}\right) _{i\in J}\right) \in \mathbb{R}%
_{+}\times \mathbb{R}_{+}^{J}$\ such that $\dsum\nolimits_{i\in
J}x_{i}^{\ast }=x^{\ast },$ $\alpha +\sum\nolimits_{i\in J}\varepsilon
_{i}=\varepsilon +\eta $ and (\ref{2}) holds.

Thus, if (\ref{1}) holds, then, for any $\varepsilon \geq 0$ we have $%
M^{\varepsilon }f(x^{\ast })\subset N^{\varepsilon }f(x^{\ast }),$ where the
multifunction $N^{\varepsilon }f:X^{\ast }\rightrightarrows X$ is defined,
for each $x^{\ast }\in X^{\ast },$ by%
\begin{equation}
N^{\varepsilon }f(x^{\ast })=\bigcap_{\eta >0}\dbigcup\limits_{\QATOP{J\in
\mathcal{F}\left( I\right) }{\hfill }}\dbigcup\limits_{\QATOP{\left(
x_{i}^{\ast }\right) _{i\in J}\in \left( X^{\ast }\right) ^{J}}{%
\sum\nolimits_{i\in J}x_{i}^{\ast }=x^{\ast }}}{{\bigcup\limits_{\QATOP{%
\left( \alpha ,\left( \varepsilon _{i}\right) _{i\in J}\right) \in \mathbb{R}%
_{+}\times \mathbb{R}_{+}^{J}}{\alpha +\sum\nolimits_{i\in J}\varepsilon
_{i}=\varepsilon +\eta }}}}\left( T_{f}^{\alpha }\left( J\right) \bigcap
\left( \bigcap_{i\in J}M^{\varepsilon _{i}}f(x_{i}^{\ast })\right) \right) .
\label{3}
\end{equation}

Since $M^{\varepsilon }f(x^{\ast })=\emptyset $ when ${{f^{\ast }}}\left(
x^{\ast }\right) \notin \mathbb{R},$ we can state:

\begin{lemma}
\label{Lemma1}If has $(\mathrm{RP}_{x^*}) $ zero duality gap at $x^{\ast
}\in X^{\ast },$ then
\begin{equation*}
M^{\varepsilon }f(x^{\ast })\subset N^{\varepsilon }f(x^{\ast }),\forall
\varepsilon \geq 0.
\end{equation*}
\end{lemma}

It turns out that the reverse inclusion always holds. Let us check this. Let
$\varepsilon \geq 0$ and $x\notin M^{\varepsilon }f(x^{\ast }).$

If $f\left( x\right) =+\infty ,$ then, by definition of $T_{f}^{\alpha
}\left( J\right) ,$ we have $T_{f}^{\alpha }\left( J\right) =\emptyset $ for
any $\left( \alpha ,J\right) \in \mathbb{R}_{+}\times \mathcal{F}\left(
I\right) .$ Consequently, $x\notin N^{\varepsilon }f(x^{\ast })=\emptyset .$

Assume now that $f\left( x\right) \in \mathbb{R}.$ Since $x\notin
M^{\varepsilon }f(x^{\ast }),$ there exists $\eta >0$ such that%
\begin{equation}
{f}\left( x\right) +{{f^{\ast }}}\left( x^{\ast }\right) -\left\langle
x^{\ast },x\right\rangle >\varepsilon +\eta .  \label{4}
\end{equation}%
Let us suppose now that $x\in N^{\varepsilon }f(x^{\ast }).$ Then, there
exist $J\in \mathcal{F}\left( I\right) ,$ $\left( x_{i}^{\ast }\right)
_{i\in J}\in \left( X^{\ast }\right) ^{J},$ and $\left( \alpha ,\left(
\varepsilon _{i}\right) _{i\in J}\right) \in \mathbb{R}_{+}\times \mathbb{R}%
_{+}^{J}$\ such that $\dsum\nolimits_{i\in J}x_{i}^{\ast }=x^{\ast },$ $%
\alpha +\sum\nolimits_{i\in J}\varepsilon _{i}=\varepsilon +\eta ,$ and $%
x\in T_{f}^{\alpha }\left( J\right) \cap \left( \bigcap_{i\in
J}M^{\varepsilon _{i}}f(x_{i}^{\ast })\right) .$ By definition of $\varphi $
we thus have ${{f^{\ast }}}\left( x^{\ast }\right) \leq \varphi \left(
x^{\ast }\right) \leq \dsum\nolimits_{i\in J}{{f_{i}^{\ast }}}\left(
x_{i}^{\ast }\right) $ and, so,%
\begin{equation*}
\begin{array}{ll}
{f}\left( x\right) +{{f^{\ast }}}\left( x^{\ast }\right) -\left\langle
x^{\ast },x\right\rangle & \leq \left[ f\left( x\right) -\sum\limits_{i\in J}%
{{f_{i}}}\left( x\right) \right] +\dsum\limits_{i\in J}\left[ {{f_{i}}}%
\left( x\right) +{{f_{i}^{\ast }}}\left( x_{i}^{\ast }\right) -\left\langle
x_{i}^{\ast },x\right\rangle \right] \\
& \leq \alpha +\sum\nolimits_{i\in J}\varepsilon _{i}=\varepsilon +\eta ,%
\end{array}%
\end{equation*}%
which contradicts (\ref{4}). So, $x\notin N^{\varepsilon }f(x^{\ast })$ and
we can claim:

\begin{lemma}
\label{Lemma2} For any $x^{\ast }\in X^{\ast },$ it holds
\begin{equation*}
N^{\varepsilon }f(x^{\ast })\subset M^{\varepsilon }f(x^{\ast }),\forall
\varepsilon \geq 0.
\end{equation*}
\end{lemma}

By Lemmas \ref{Lemma1} and \ref{Lemma2} we have

\begin{lemma}
\label{Lemma3}Let $x^{\ast }\in X^{\ast }$ be such that ${{f^{\ast }}}\left(
x^{\ast }\right) =\varphi \left( x^{\ast }\right) .$ Then,%
\begin{equation*}
M^{\varepsilon }f(x^{\ast })=N^{\varepsilon }f(x^{\ast }),\forall
\varepsilon \geq 0.
\end{equation*}
\end{lemma}

It turns out again that the reverse implication always holds. In fact, we
can prove a little bit more:

\begin{lemma}
\label{Lemma4}Let $x^{\ast }\in X^{\ast }$ and assume that there exists $%
\overline{\varepsilon }>0$ such that%
\begin{equation*}
M^{\varepsilon }f(x^{\ast })\subset N^{\varepsilon }f(x^{\ast }),\forall
\varepsilon \in \left] 0,\overline{\varepsilon }\right[ .
\end{equation*}%
Then ${{f^{\ast }}}\left( x^{\ast }\right) =\varphi \left( x^{\ast }\right)
. $
\end{lemma}

\textit{Proof. }We have just to check that $\varphi \left( x^{\ast }\right)
\leq {{f^{\ast }}}\left( x^{\ast }\right) .$ This is obvious if ${{f^{\ast }}%
}\left( x^{\ast }\right) =+\infty .$ Assume now that ${{f^{\ast }}}\left(
x^{\ast }\right) \in \mathbb{R}.$ Let us assume that $\varphi \left( x^{\ast
}\right) >{{f^{\ast }}}\left( x^{\ast }\right) .$

There exists $\varepsilon \in \left] 0,\overline{\varepsilon }\right[ $ such
that
\begin{equation}
\varphi \left( x^{\ast }\right) >{{f^{\ast }}}\left( x^{\ast }\right)
+3\varepsilon .  \label{5}
\end{equation}

Let us pick $x\in M^{\varepsilon }f(x^{\ast }),$ which is non-empty since ${{%
f^{\ast }}}\left( x^{\ast }\right) \in \mathbb{R}.$ By hypothesis $x\in
N^{\varepsilon }f(x^{\ast })$ and, by (\ref{3}), with $\eta =\varepsilon ,$
there exist $J\in \mathcal{F}\left( I\right) ,$ $\left( x_{i}^{\ast }\right)
_{i\in J}\in \left( X^{\ast }\right) ^{J},$ and $\left( \alpha ,\left(
\varepsilon _{i}\right) _{i\in J}\right) \in \mathbb{R}_{+}\times \mathbb{R}%
_{+}^{J}$\ such that $\dsum\nolimits_{i\in J}x_{i}^{\ast }=x^{\ast },$ $%
\alpha +\sum\nolimits_{i\in J}\varepsilon _{i}=2\varepsilon ,$ $f\left(
x\right) -\sum\nolimits_{i\in J}{{f_{i}}}\left( x\right) \leq \alpha ,$\ and
${{f_{i}}}\left( x\right) +{{f_{i}^{\ast }}}\left( x_{i}^{\ast }\right) \leq
\left\langle x_{i}^{\ast },x\right\rangle +\varepsilon _{i},$ for all $i\in
J.$ We thus have
\begin{equation*}
\begin{array}{ll}
-{{f^{\ast }}}\left( x^{\ast }\right) & \leq f\left( x\right) -\left\langle
x^{\ast },x\right\rangle \\
& =\left[ f\left( x\right) -\sum\limits_{i\in J}{{f_{i}}}\left( x\right) %
\right] +\dsum\limits_{i\in J}\left[ {{f_{i}}}\left( x\right) +{{f_{i}^{\ast
}}}\left( x_{i}^{\ast }\right) -\left\langle x_{i}^{\ast },x\right\rangle %
\right] -\dsum\limits_{i\in J}{{f_{i}^{\ast }}}\left( x_{i}^{\ast }\right)
\\
& \leq \alpha +\sum\limits_{i\in J}\varepsilon _{i}-\dsum\limits_{i\in J}{{%
f_{i}^{\ast }}}\left( x_{i}^{\ast }\right) \\
& \leq 2\varepsilon -\varphi \left( x^{\ast }\right) ,%
\end{array}%
\end{equation*}%
which contradicts (\ref{5}). So, $\varphi \left( x^{\ast }\right) \leq {{%
f^{\ast }}}\left( x^{\ast }\right)$, which together with the weak duality
shows that $\varphi \left( x^{\ast }\right) = {{f^{\ast }}}\left( x^{\ast
}\right)$ and we are done. \qquad \qed

We now state the main result of this section.

\begin{theorem}[Zero duality gap]
\label{Th1}Let $\left( f_{i}\right) _{i\in I}$ be a family of proper
functions with $f=\sum\nolimits_{i\in I}^{R}{{f_{i}}}$ proper, and let ${x}%
^{\ast }\in X^{\ast }.$ The next statements are equivalent:\newline
$(i)$ $(\mathrm{RP}_{x^*}) $ has zero duality gap, \newline
$(ii)$ $M^{\varepsilon }f({x}^{\ast })=N^{\varepsilon }f({x }^{\ast
}),\forall \varepsilon \geq 0, $\newline
$(iii)$ There exists $\overline{\varepsilon }>0$ such that%
\begin{equation*}
M^{\varepsilon }f({x}^{\ast })=N^{\varepsilon }f({x}^{\ast }),\forall
\varepsilon \in \left] 0,\overline{\varepsilon }\right[ ,
\end{equation*}%
$(iv)$ There exists $\overline{\varepsilon }>0$ such that%
\begin{equation*}
M^{\varepsilon }f({x}^{\ast })\subset N^{\varepsilon }f({x}^{\ast }),\forall
\varepsilon \in \left] 0,\overline{\varepsilon }\right[.
\end{equation*}
\end{theorem}

\textit{Proof.} Lemma \ref{Lemma3} says that $\left[ (i)\Longrightarrow (ii)%
\right] ,$ while $\left[ (ii)\Longrightarrow (iii)\right] $ and $\left[
(iii)\Longrightarrow (iv)\right] $ are obvious. Finally, $\left[
(vi)\Longrightarrow (i)\right] $ is Lemma \ref{Lemma4}.\qquad \qed

We now characterize stable zero duality gap for the robust sum problem. To
this end, let us introduce $\Pi ^{\varepsilon }f:=\left( N^{\varepsilon
}f\right) ^{-1},$ i.e., the inverse multifunction of $N^{\varepsilon }f.$
One has $\Pi ^{\varepsilon }\dot{f}:X\rightrightarrows X^{\ast }$ and, for
any $\left( x^{\ast },x\right) \in X^{\ast }\times X,$%
\begin{equation*}
x^{\ast }\in \Pi ^{\varepsilon }f\left( x\right) \Longleftrightarrow x\in
N^{\varepsilon }f\left( x^{\ast }\right) .
\end{equation*}%
The next explicit formula holds:

\begin{lemma}
\label{Lemma5}For any $\left( x,\varepsilon \right) \in X\times \mathbb{R}%
_{+}$ we have%
\begin{equation}
\Pi ^{\varepsilon }f\left( x\right) =\bigcap_{\eta >0}\dbigcup\limits_{%
\QATOP{0\leq \alpha \leq \varepsilon +\eta }{{}}}\dbigcup\limits_{\QATOP{%
J\in S_{f}^{\alpha }\left( x\right) }{\hfill }}{{\bigcup\limits_{\QATOP{%
\scriptstyle\left( \varepsilon _{i}\right) _{i\in J}\in \mathbb{R}_{+}^{J}}{%
\scriptstyle\sum\nolimits_{i\in J}\varepsilon _{i}=\varepsilon +\eta -\alpha
}}}}\dsum\limits_{i\in J}\partial ^{\varepsilon _{i}}f_{i}\left( x\right) .
\label{6}
\end{equation}
\end{lemma}

\textit{Proof. }By (\ref{3}) we have $x^{\ast }\in \Pi ^{\varepsilon
}f\left( x\right) $ if and only if for any $\eta >0$ there exist $J\in
\mathcal{F}\left( I\right) ,$ $\left( x_{i}^{\ast }\right) _{i\in J}\in
\left( X^{\ast }\right) ^{J},$ and $\left( \alpha ,\left( \varepsilon
_{i}\right) _{i\in J}\right) \in \mathbb{R}_{+}\times \mathbb{R}_{+}^{J}$\
such that $\dsum_{i\in J}x_{i}^{\ast }=x^{\ast },$ $\alpha
+\sum\nolimits_{i\in J}\varepsilon _{i}=\varepsilon +\eta ,$ $x\in
T_{f}^{\alpha }\left( J\right) $ (i.e., $J\in S_{f}^{\alpha }\left( x\right)
$), and $x\in \bigcap_{i\in J}M^{\varepsilon _{i}}f_{i}(x_{i}^{\ast }) $
(i.e., $x_{i}^{\ast }\in \partial ^{\varepsilon _{i}}f_{i}\left( x\right) $
for all $i\in J$). This exactly means that $x^{\ast }$ belongs to the set in
the right hand side of (\ref{6}). \qquad \qed

\begin{lemma}
\label{Lemma6}For any $\left( x,\varepsilon \right) \in X\times \mathbb{R}%
_{+}$ one has%
\begin{equation*}
\Pi ^{\varepsilon }f\left( x\right) \subset \partial ^{\varepsilon }f\left(
x\right) .
\end{equation*}
\end{lemma}

\textit{Proof. }Let $x^{\ast }\in \Pi ^{\varepsilon }f\left( x\right) .$ We
have $x\in N^{\varepsilon }f\left( x^{\ast }\right) $ and, by Lemma \ref%
{Lemma2}, $x\in M^{\varepsilon }f(x^{\ast }),$ that means $x^{\ast }\in
\partial ^{\varepsilon }f\left( x\right) .$ \qquad \qed

We now characterize the stable zero duality gap for the robust sum problem.

\begin{theorem}[Stable zero duality gap]
\label{Th2} Let $\left( f_{i}\right) _{i\in I}$ be a family of proper
functions with $f=\sum\nolimits_{i\in I}^{R}{{f_{i}}\left( x\right) }$
proper. The next statements are equivalent:\newline
$(i)$ ${{f^{\ast }}}\left( x^{\ast }\right) =\varphi \left( x^{\ast }\right)
,\forall x^{\ast }\in X^{\ast }, $\newline
$(ii)$ $\partial ^{\varepsilon }f\left( x\right) =\Pi ^{\varepsilon }f\left(
x\right) ,\forall \left( x,\varepsilon \right) \in X\times \mathbb{R}_{+}, $%
\newline
$(iii)$ There exists $\overline{\varepsilon }>0$ such that%
\begin{equation*}
\partial ^{\varepsilon }f\left( x\right) =\Pi ^{\varepsilon }f\left(
x\right) ,\forall \left( x,\varepsilon \right) \in X\times \left] 0,%
\overline{\varepsilon }\right[,
\end{equation*}%
$(iv)$ There exists $\overline{\varepsilon }>0$ such that%
\begin{equation*}
\partial ^{\varepsilon }f\left( x\right) \subset \Pi ^{\varepsilon }f\left(
x\right) ,\forall \left( x,\varepsilon \right) \in X\times \left] 0,%
\overline{\varepsilon }\right[.
\end{equation*}%
\qquad
\end{theorem}

\textit{Proof. } $\left[ (i)\Longrightarrow (ii)\right] $ Let $\left(
x,\varepsilon \right) \in X\times \mathbb{R}_{+}.$ We know that $x^{\ast
}\in \partial ^{\varepsilon }f\left( x\right) $ if and only if $x\in
M^{\varepsilon }f(x^{\ast }).$ By Theorem \ref{Th1}, $M^{\varepsilon
}f(x^{\ast })=N^{\varepsilon }f(x^{\ast }).$ So,%
\begin{equation*}
x^{\ast }\in \partial ^{\varepsilon }f\left( x\right) \Longleftrightarrow
x\in N^{\varepsilon }f(x^{\ast })\Longleftrightarrow x^{\ast }\in \Pi
^{\varepsilon }f\left( x\right)
\end{equation*}%
and $(ii)$ holds.

$\left[ (ii)\Longrightarrow (iii)\right] $ and $\left[ (iii)\Longrightarrow
(iv)\right] $ are obvious.

$\left[ (vi)\Longrightarrow (i)\right] $ Let $\left( x^{\ast },\varepsilon
\right) \in X^{\ast }\times \left] 0,\overline{\varepsilon }\right[ $ and $%
x\in M^{\varepsilon }f(x^{\ast }).$ We have $x^{\ast }\in \partial
^{\varepsilon }f\left( x\right) $ and, by Theorem \ref{Th1}, $x^{\ast }\in
\Pi ^{\varepsilon }f\left( x\right) ,$ that means $x\in N^{\varepsilon
}f(x^{\ast }).$ So, $M^{\varepsilon }f(x^{\ast })\subset N^{\varepsilon
}f(x^{\ast })$ for any $\varepsilon \in \left] 0,\overline{\varepsilon }%
\right[ ,$ and, again by Theorem \ref{Th1}, ${{f^{\ast }}}\left( x^{\ast
}\right) =\varphi \left( x^{\ast }\right) .$ \qquad \qed

\section{Strong duality\protect\bigskip}

\begin{definition}
We say that the robust sum problem $(\mathrm{RP}_{x^*}) $ has a strong zero
duality gap at a given $x^{\ast }\in X^{\ast }$ if there exist $J\in
\mathcal{F}\left( I\right) $ and $\left( x_{i}^{\ast }\right) _{i\in J}\in
\left( X^{\ast }\right) ^{J}$ such that $x^{\ast }=\dsum\limits_{i\in
J}x_{i}^{\ast }$ and

\begin{equation}
\inf (\mathrm{RP}_{x^*}) = - f^{\ast }\left( x^{\ast }\right) = -
\dsum\limits_{i\in J}{{f_{i}^{\ast }}}\left( x_{i}^{\ast }\right) = \sup (%
\mathrm{RD}_{x^*}).
\end{equation}
If the above condition holds at each $x^{\ast }\in X^{\ast }$ we will say
that $(\mathrm{RP}_{x^*}) $ has a stable strong zero duality gap.
\end{definition}

To characterize the strong zero duality gap of the robust sum problem $(%
\mathrm{RP}_{x^{\ast }})$, let us fix some notation first. Given $x^{\ast
}\in X^{\ast }$, $\varepsilon \geq 0$, $J\in \mathcal{F}(I)$, $(x_{i}^{\ast
})_{\imath \in J}\in (X^{\ast })^{J}$, define
\begin{equation*}
B_{(J,(x_{i}^{\ast })_{i\in J})}^{\varepsilon }f(x^{\ast }):=\left\{
\begin{array}{ll}
{{\bigcup\limits_{\QATOP{\left( \alpha ,\left( \varepsilon _{i}\right)
_{i\in J}\right) \in \mathbb{R}_{+}\times \mathbb{R}_{+}^{J}}{\alpha
+\sum\nolimits_{i\in J}\varepsilon _{i}=\varepsilon }}}}T_{f}^{\alpha
}\left( J\right) \bigcap \left( \bigcap_{i\in J}M^{\varepsilon
_{i}}f(x_{i}^{\ast })\right) , & \ \ \text{\textrm{if}}\sum\limits_{i\in
J}x_{i}^{\ast }=x^{\ast }, \\
\emptyset , & \ \text{else.}%
\end{array}%
\right.
\end{equation*}

\begin{theorem}[Strong zero duality gap]
\label{Th3}Let $\left( f_{i}\right) _{i\in I}$ be a family of proper
functions with $f=\sum\nolimits_{i\in I}^{R}{{f_{i}}}$ proper, and let ${x}%
^{\ast }\in X^{\ast }.$ The next statements are equivalent:\newline
$(i)$ The robust sum problem $(\mathrm{RP}_{x^*}) $ has a strong zero
duality gap, \newline
$(ii)$ $\exists J \in \mathcal{F}(I)$, $\exists (x^*_i)_{i\in J} \in (X^\ast
)^J$: $M^\varepsilon f(x^*) = B^{\varepsilon }_{(J, (x^*_i)_{i \in
J})}f(x^*) $, $\forall \varepsilon \geq 0, $ \newline
$(iii)$ There exist $\overline{\varepsilon }>0$, $J \in \mathcal{F}(I)$, $%
(x^*_i)_{i\in J} \in (X^\ast)^J $ such that%
\begin{equation}  \label{new2}
M^{\varepsilon }f({x}^{\ast })= B^{\varepsilon }_{(J, (x^*_i)_{i \in
J})}f(x^*) , \ \ \forall \varepsilon \in \left] 0,\overline{\varepsilon }%
\right[,
\end{equation}
\end{theorem}

\textit{Proof.} $\left[ (i)\Longrightarrow (ii)\right] $ By the very
definition of $B_{(J,(x_{i}^{\ast })_{i\in J})}^{\varepsilon }f(x^{\ast })$, %
\eqref{3}, and Lemma \ref{Lemma2} we have
\begin{equation*}
B_{(J,(x_{i}^{\ast })_{i\in J})}^{\varepsilon }f(x^{\ast })\subset
N^{\varepsilon }f(x^{\ast })\subset M^{\varepsilon }f(x^{\ast }).
\end{equation*}%
Let $x\in M^{\varepsilon }f(x^{\ast })$. By $(i)$ there exist $J\in \mathcal{%
F}(I)$, $(x_{i}^{\ast })_{i\in J}\in (X^{\ast })^{J}$ such that $\sum_{i\in
J}x_{i}^{\ast }=x^{\ast }$ and
\begin{equation*}
\sum\limits_{i\in J}f_{i}^{\ast }(x_{i}^{\ast })=f^{\ast }(x^{\ast })\leq
\langle x^{\ast },x\rangle -f(x)+\varepsilon .
\end{equation*}%
Consequently,
\begin{equation*}
\sum\limits_{i\in J}\Big[f_{i}^{\ast }(x_{i}^{\ast })+f_{i}(x)-\langle
x_{i}^{\ast },x\rangle \Big]+\Big[f(x)-\sum\limits_{i\in J}f_{i}(x)\Big]\leq
\varepsilon .
\end{equation*}%
Since all the above brackets are non negative, there exist $(\alpha
,(\varepsilon _{i})_{i})\in \mathbb{R}_{+}\times \mathbb{R}_{+}^{J}$ such
that $f(x)-\sum_{i\in J}f_{i}(x)\leq \alpha $, that means $x\in
T_{f}^{\alpha }(J)$, $\alpha +\sum_{i\in J}\varepsilon _{i}=\varepsilon $,
and for each $i\in J$,
\begin{equation*}
f_{i}^{\ast }(x_{i}^{\ast })+f_{i}(x)-\langle x_{i}^{\ast },x\rangle \leq
\varepsilon _{i},
\end{equation*}%
that means $x\in \cap _{i\in J}M^{\varepsilon _{i}}f_{i}(x_{i}^{\ast })$. So
$x\in B_{(J,(x_{i}^{\ast })_{i\in J})}^{\varepsilon }f(x^{\ast })$ and $(ii)$
holds.

$\left[ (ii)\Longrightarrow (iii)\right] $ is obvious.

$\left[ (iii)\Longrightarrow (i)\right] $ Assume that $(iii)$ holds. So,
there exist $\overline{\varepsilon }>0$, $J \in \mathcal{F}(I)$, $(x^*_i)_{i
\in J} \in (X^\ast)^J $ such that \eqref{new2} holds. Let us first prove
that $\sum_{i \in J} f^\ast_i (x^*_i) \leq f^\ast (x^*)$. Assume the
contrary, i.e., there exists $\varepsilon > 0$, that we can choose $%
\varepsilon < \overline{\varepsilon} $, such that
\begin{equation}  \label{new1}
f^\ast (x^*) + \varepsilon < \sum_{i \in J} f^\ast_i(x^*_i).
\end{equation}
We have $f^\ast (x^*) \in \mathbb{R}$. Picking $x \in M^\varepsilon f(x^*)$
which is non-empty, we have $x \in B^{\varepsilon }_{(J, (x^*_i)_{i \in
J})}f(x^*) $ and hence, there exist $(\alpha, (\varepsilon_i)_i ) \in
\mathbb{R}_+\times \mathbb{R}_+^J$ such that $\alpha + \sum_{i \in J}
\varepsilon_i = \varepsilon $, $\sum_{i \in J} x^*_i = x^*$ and $x \in
T^\alpha_f(J) \cap\left( \cap_{i \in J} M^{\varepsilon_i} {f_i}
(x^*_i)\right)$. Then
\begin{eqnarray*}
\sum_{i \in J} f^\ast_i(x^*_i) &\leq& \sum_{i\in J}\left[ \langle x^*_i,
x\rangle - f_i(x) + \varepsilon_i\right] = \langle x^*, x\rangle -
\sum_{i\in J}f_i(x) + \sum_{i\in J}\varepsilon_i \\
&\leq& \langle x^*, x\rangle - f(x) + \alpha + \sum_{i \in J} \varepsilon_i
= \langle x^*, x\rangle - f(x) + \varepsilon \\
&\leq& f^\ast(x^*) + \varepsilon,
\end{eqnarray*}
which contradicts \eqref{new1}. We then have
\begin{equation*}
\varphi (x^*) \leq \sum_{i \in J}f_i^\ast (x^*_i) \leq f^\ast (x^*) \leq
\varphi(x^*).
\end{equation*}
So, $\varphi (x^*) = \sum_{i\in J} f_i^\ast (x^*_i) = f^\ast(x^*)$ with $%
\sum_{i \in J} x^*_i = x^*$, that means that $(i)$ holds. \qquad \qed

In order to characterize the stable strong zero duality gap for the robust
sum problem $(\mathrm{RP}_{x^*}) $, let us introduce, for each $\varepsilon
\geq 0$, the set-valued mapping $N_{s}^{\varepsilon }f:X^{\ast
}\rightrightarrows X$ defined by
\begin{equation*}
N_{s}^{\varepsilon }f(x^{\ast }) := \dbigcup\limits_{\QATOP{J\in \mathcal{F}%
\left( I\right), \left( x_{i}^{\ast }\right) _{i\in J}\in \left( X^{\ast
}\right) ^{J}}{\sum\nolimits_{i\in J}x_{i}^{\ast }=x^{\ast }}}
B^{\varepsilon }_{(J, (x^*_i)_{i \in J})}f(x^*) , \ \forall x^* \in X^\ast ,
\end{equation*}%
and its inverse $\Pi _{s}^{\varepsilon }f : X \rightrightarrows X^\ast $.
For each $(x, x^*) \in X\times X^\ast$ one has
\begin{equation*}
{x}^{\ast }\in \Pi _{s}^{\varepsilon }f\left( x\right) \Longleftrightarrow
x\in N_{s}^{\varepsilon }f\left( {x}^{\ast }\right) .
\end{equation*}
More explicitly one has straightforwardly, for each $x \in X$,
\begin{equation*}
\Pi _{s}^{\varepsilon }f\left( x\right) =\dbigcup\limits_{\QATOP{0\leq
\alpha \leq \varepsilon }{{}}}\dbigcup\limits_{\QATOP{J\in S_{f}^{\alpha
}\left( x\right) }{\hfill }}{{\bigcup\limits_{\QATOP{\scriptstyle\left(
\varepsilon _{i}\right) _{i\in J}\in \mathbb{R}_{+}^{J}}{\scriptstyle%
\sum\nolimits_{i\in J}\varepsilon _{i}=\varepsilon -\alpha }}}}%
\dsum\limits_{i\in J}\partial ^{\varepsilon _{i}}f_{i}\left( x\right),
\end{equation*}%
where $S_{f}^{\alpha}\left( x\right) = \{ J \in \mathcal{F}(I) \, :\,
\sum_{i\in J} f_i(x) +\alpha \geq f(x) \in \mathbb{R}\}$ as in Section 4. We
have
\begin{equation*}
N_{s}^{\varepsilon }f(x^{\ast }) \subset (N^\varepsilon f) (x^*) \subset
(M^\varepsilon f)(x^*),
\end{equation*}
and, passing to the inverse multivalued mappings,
\begin{equation}  \label{eq5.4}
(\Pi^\varepsilon_sf)(x) \subset (\Pi^\varepsilon f)(x) \subset \partial
^\varepsilon f (x) , \forall x \in X, \forall \varepsilon \geq 0.
\end{equation}

\begin{theorem}[Stable strong zero duality gap]
\label{Th4} Let $\left( f_{i}\right) _{i\in I}$ be a family of proper
functions with $f=\sum\nolimits_{i\in I}^{R}{{f_{i}}}$ proper. The next
statements are equivalent:\newline
$(i)$ The robust sum problem $(\mathrm{RP}_{x^*}) $ has stable strong zero
duality gap, \newline
$(ii)$ $\partial ^{\varepsilon }f\left( x\right) =\Pi _{s}^{\varepsilon
}f\left( x\right) ,\forall \left( x,\varepsilon \right) \in X\times \mathbb{R%
}_{+}$, \newline
$(iii)$ $\exists \overline{\varepsilon }>0$$:$ $\partial ^{\varepsilon
}f\left( x\right) =\Pi _{s}^{\varepsilon }f\left( x\right) , \forall \left(
x,\varepsilon \right) \in X\times \left[ 0,\overline{\varepsilon }\right].$
\end{theorem}

\textit{Proof.} $\left[ (i)\Longrightarrow (ii)\right] $ We only have to
prove the inclusion \textquotedblleft $\subset $" in $(ii)$. So, let $%
x^{\ast }\in \partial ^{\varepsilon }f(x)$. By $(i)$, there exist $J\in
\mathcal{F}(I)$, $(x_{i}^{\ast })_{i\in J}\in (X^{\ast })^{J}$ such that $%
\sum_{i\in J}x_{i}^{\ast }=x^{\ast }$ and
\begin{equation*}
\sum_{i\in J}f_{i}^{\ast }(x_{i}^{\ast })=f^{\ast }(x^{\ast })\leq \langle
x^{\ast },x\rangle -f(x)+\varepsilon .
\end{equation*}%
Consequently,
\begin{equation*}
\sum_{i\in J}\Big[f_{i}^{\ast }(x_{i}^{\ast })+f_{i}(x)-\langle x_{i}^{\ast
},x\rangle \Big]+\Big[f(x)-\sum_{i\in J}f_{i}^{\ast }(x_{i}^{\ast })\Big]%
\leq \varepsilon ,
\end{equation*}%
and there exist $((\varepsilon _{i})_{i\in J},\alpha )\in \mathbb{R}%
^{J}\times \mathbb{R}$ such that $f_{i}^{\ast }(x_{i}^{\ast
})+f_{i}(x)-\langle x_{i}^{\ast },x\rangle \leq \varepsilon _{i}$ for each $%
i\in J$, $f(x)-\sum_{i\in J}f_{i}(x)\leq \alpha $, and $\alpha +\sum_{i\in
J}\varepsilon _{i}=\varepsilon $. We thus have $x^{\ast }=\sum_{i\in
J}x_{i}^{\ast }\in \sum_{i\in J}\partial ^{\varepsilon }f_{i}(x)$ with $%
0\leq \alpha \leq \varepsilon $, $J\in S_{f}^{\alpha }(x)$, and $\sum_{i\in
J}\varepsilon _{i}=\varepsilon -\alpha $, that means $x^{\ast }\in (\Pi
_{s}^{\varepsilon }f)(x)$, and $(ii)$ holds.

$\left[ (ii)\Longrightarrow (iii)\right] $ is obvious.

$\left[ (iii)\Longrightarrow (i)\right] $ Let $x^{\ast }\in X^{\ast }$. If $%
f^{\ast }(x^{\ast })=+\infty $ then $\varphi (x^{\ast })=+\infty $ and $f$
has obviously a strong zero duality gap at $x^{\ast }$. Since $\limfunc{dom}%
f\not=\emptyset $ we have $f^{\ast }(x^{\ast })\neq -\infty $ and it remains
to consider the case $f^{\ast }(x^{\ast })\in \mathbb{R}$. Pick $x\in M^{%
\overline{\varepsilon }}f(x^{\ast })$ which is non-empty, and set $%
\varepsilon :=f^{\ast }(x^{\ast })+f(x)-\langle x^{\ast },x\rangle $. One
has $\varepsilon \in \lbrack 0,\overline{\varepsilon }]$, $x^{\ast }\in
\partial ^{\varepsilon }f(x)$ and, by $(iii),$ there exist $\alpha \in
\lbrack 0,\varepsilon ]$, $J\in S_{f}^{\alpha }(x)$, $(x_{i}^{\ast })_{i\in
J}\in (X^{\ast })^{J}$, $(\varepsilon _{i})_{i\in J}\in \mathbb{R}_{+}^{J}$
such that $\alpha +\sum_{i\in J}\varepsilon _{i}=\varepsilon $, $\sum_{i\in
J}x_{i}^{\ast }=x^{\ast }$, and $x_{i}^{\ast }\in \partial ^{\varepsilon
_{i}}f_{i}(x)$ for each $i\in J$. We thus have
\begin{eqnarray*}
\varphi (x^{\ast }) &\leq &\sum_{i\in J}f_{i}^{\ast }(x_{i}^{\ast })\leq
\sum_{i\in J}\Big[\langle x_{i}^{\ast },x\rangle -f_{i}(x)+\varepsilon _{i}%
\Big] \\
&=&\langle x^{\ast },x\rangle -\sum_{i\in J}f_{i}(x)+\sum_{i\in
J}\varepsilon _{i} \\
&\leq &\langle x^{\ast },x\rangle -f(x)+\alpha +\sum_{i\in J}\varepsilon _{i}
\\
&=&\langle x^{\ast },x\rangle -f(x)+\varepsilon =f^{\ast }(x^{\ast })\leq
\varphi (x^{\ast }).
\end{eqnarray*}%
Consequently, $f^{\ast }(x^{\ast })=\sum_{i\in J}f_{i}^{\ast }(x_{i}^{\ast
}) $ with $J\in \mathcal{F}(I)$ and $\sum_{i\in J}x_{i}^{\ast }=x^{\ast }$,
that means $f^{\ast }$ has strong zero duality gap at $x^{\ast }$ and we are
done. \qquad \qed

\section{Duality for the robust sum of closed convex functions}

Denote by $\limfunc{co}A$ the convex hull of $A\subset X^{\ast }\times
\mathbb{R}$, by $\overline{A}$ its closure w.r.t. the $w^{\ast }-$topology
and by $\overline{\limfunc{co}}A$ its $w^{\ast }-$closed convex hull. We
also denote by $\Gamma \left( X\right) $ the set of all proper convex lsc
functions on $X.$ In this section we assume that
\begin{equation}
\left( f_{i}\right) _{i\in I}\subset \Gamma \left( X\right) \text{ and }
\limfunc{dom } f \not= \emptyset  \label{7}
\end{equation}
(recall that $f=\sum\nolimits_{i\in I}^{R}{{f_{i}}}$). We thus have $f \in
\Gamma(X)$.

Let us introduce the set
\begin{equation*}
\mathcal{A}:=\dbigcup\limits_{\QATOP{J\in \mathcal{F}\left( I\right) }{%
\hfill }}\sum\nolimits_{i\in J}\limfunc{epi}f_{i}^{\ast },
\end{equation*}%
which is related with the function%
\begin{equation*}
\varphi \left( x^{\ast }\right) :=\inf_{J\in \mathcal{F}\left( I\right)
}\left\{ \dsum\limits_{i\in J}{{f_{i}^{\ast }}}\left( x^{\ast }\right)
:\left( x_{i}^{\ast }\right) _{i\in J}\in \left( X^{\ast }\right)
^{J},\dsum\limits_{i\in J}x_{i}^{\ast }=x^{\ast }\right\} ,\forall x^{\ast
}\in X^{\ast },
\end{equation*}%
by the (easily checkable) double inclusion%
\begin{equation}
\limfunc{epi}\nolimits_{s}\varphi \subset \mathcal{A\subset }\limfunc{epi}%
\varphi .  \label{5.1}
\end{equation}%
Thus,
\begin{equation}  \label{6.2b}
\overline{\limfunc{co}}\mathcal{A} \ = \ \overline{\limfunc{co}}\limfunc{epi}%
\varphi .
\end{equation}

\begin{lemma}
\label{Lemma7} Assume that (\ref{7}) holds. Then $\varphi ^{\ast }=f$ and $%
\limfunc{epi}f^{\ast }=\overline{\limfunc{co}}\mathcal{A}.$
\end{lemma}

\textit{Proof.} We have $\varphi =\inf\limits_{\QATOP{J\in \mathcal{F}\left(
I\right) }{\hfill }}\left( \square_{i\in J}f_{i}^{\ast }\right) ,$ where
\begin{equation*}
\square_{i\in J}f_{i}^{\ast }\left( x^{\ast }\right) :=\inf \left\{
\dsum\limits_{i\in J}f_{i}^{\ast }\left( x_{i}^{\ast }\right)
:\dsum\limits_{i\in J}x_{i}^{\ast }=x^{\ast }\right\} ,\forall x^{\ast }\in
X^{\ast },
\end{equation*}%
is the infimal convolution of the finite family of functions $\left\{
f_{i}^{\ast },i\in J\right\} .$ So,%
\begin{equation*}
\varphi ^{\ast }=\sup\limits_{J\in \mathcal{F}\left( I\right) }\left(
\square _{i\in J}f_{i}^{\ast }\right) ^{\ast }=\sup\limits_{J\in \mathcal{F}%
\left( I\right) }\dsum\limits_{i\in J}f_{i}^{\ast \ast }=\sup\limits_{J\in
\mathcal{F}\left( I\right) }\dsum\limits_{i\in J}f_{i}=f.
\end{equation*}
For the second statement, one has $f^{\ast }=\varphi ^{\ast \ast }$ and,
since $f^{\ast }$ is proper, $\limfunc{epi} f^\ast = \limfunc{epi}\varphi
^{\ast \ast }=\overline{\limfunc{co}}\limfunc{epi}\varphi =\overline{%
\limfunc{co}}\mathcal{A}$ (the last equality follows from \eqref{6.2b}).
\qquad \qed

To go further let us recall the following notions (see, e.g., \cite{Bot10},
\cite{DGLM17}, and \cite{EV16}).

\begin{definition}
A subset $A\subset X^{\ast }\times \mathbb{R}$ is said to be closed
(respectively, closed convex) regarding another subset $B\subset X^{\ast
}\times \mathbb{R}$ if $B\cap \overline{A}$ $=B\cap A$ (respectively, $B\cap
\overline{\limfunc{co}}A=B\cap A$).
\end{definition}

\begin{theorem}[Strong zero duality gap under convexity]
\label{Th5}Assume that (\ref{7}) holds and let ${x}^{\ast }\in X^{\ast }.$
The next statements are equivalent:\newline
$(i)$ The robust sum problem $(\mathrm{RP}_{x^*}) $ has a strong zero
duality gap, \newline
$(ii)$ $\mathcal{A}$ is closed convex regarding $\left\{{x}^{\ast }\right\}
\times \mathbb{R}.$
\end{theorem}

\textit{Proof.} Assume that $f^{\ast }\left( {x}^{\ast }\right) =+\infty .$
By Lemma \ref{Lemma7}, we have $\left( \left\{{x}^{\ast }\right\} \times
\mathbb{R}\right) \cap \overline{\limfunc{co}}\mathcal{A}=\emptyset $ and $%
(ii)$ holds. By\ Proposition \ref{Prop1}, $\varphi \left({x}^{\ast }\right)
=+\infty $ and $(i)$ holds too. So, in this case, both statements $(i)$ and$%
\ (ii)$ hold.

We now assume that $f^{\ast }\left({x}^{\ast }\right) <+\infty .$ Since $f$
is proper we have $r:=f^{\ast }\left({x}^{\ast }\right) \in \mathbb{R}$ and,
by Lemma \ref{Lemma7}, $\left( {x}^{\ast },r\right) \in \overline{\limfunc{co%
}}\mathcal{A}.$

Assume now that $(ii)$ holds. Then $\left({x}^{\ast },r\right) \in \mathcal{A%
}$ and there exist $J\in \mathcal{F}\left( I\right) ,$ $\left( x_{i}^{\ast
}\right) _{i\in J}\in \left( X^{\ast }\right) ^{J},$ and $\left(
r_{i}\right) _{i\in J}\in \mathbb{R}^{J}$ such that $\dsum\limits_{i\in
J}x_{i}^{\ast }= {x}^{\ast },$ $\dsum\limits_{i\in J}r_{i}=r,$ and ${{%
f_{i}^{\ast }}}\left( x_{i}^{\ast }\right) \leq r_{i}$ for all $i\in J.$
Then, again by Proposition \ref{Prop1}, we have%
\begin{equation*}
\varphi \left( {x}^{\ast }\right) \leq \dsum\limits_{i\in J}{{f_{i}^{\ast }}}%
\left( x_{i}^{\ast }\right) \leq \dsum\limits_{i\in J}r_{i}= r = f^{\ast
}\left({x}^{\ast }\right) \leq \varphi \left( {x}^{\ast }\right) ,
\end{equation*}%
that means that $(i)$ holds.

To conclude the proof assume now that $(i)$ holds. Let $r\in \mathbb{R}$ be
such that $\left({x}^{\ast },r\right) \in \overline{\limfunc{co}}\mathcal{A}%
. $ By Lemma \ref{Lemma7} we have $f^{\ast }\left({x}^{\ast }\right) \leq r$
and, by $(i)$ there exist $J\in \mathcal{F}\left( I\right) $ and $\left(
x_{i}^{\ast }\right) _{i\in J}\in \left( X^{\ast }\right) ^{J}$ such that $%
\dsum\nolimits_{i\in J}x_{i}^{\ast }={x}^{\ast }$ and $f^\ast (x^*) =
\dsum\nolimits_{i\in J}{{f_{i}^{\ast }}}\left( x_{i}^{\ast }\right) \leq r.$
\ From this last inequality, there exists $\left( r_{i}\right) _{i\in J}\in
\mathbb{R}^{J}$ such that ${{f_{i}^{\ast }}}\left( x_{i}^{\ast }\right) \leq
r_{i},$ for all $i\in J,$ and $\dsum\nolimits_{i\in J}r_{i}=r.$ It follows
that%
\begin{equation*}
\left({x}^{\ast },r\right) =\dsum\limits_{i\in J}\left( x_{i}^{\ast
},r_{i}\right) \in \dsum\limits_{i\in J}\limfunc{epi}f_{i}^{\ast }\in
\mathcal{A}.\text{ \ \ \ \ \ \ \qed}
\end{equation*}

Since $\mathcal{A}$\ is closed convex if and only if it is closed convex
regarding $\left\{ x^{\ast }\right\} \times \mathbb{R}$ for all $x^{\ast
}\in X^{\ast },$ we have:

\begin{corollary}[Stable strong zero duality gap]
\label{Corol1}Assume that (\ref{7}) holds. The next statements are
equivalent:\newline
$(i)$ The robust sum problem $(\mathrm{RP}_{x^*}) $ has stable strong zero
duality gap, \newline
$(ii)$ $\mathcal{A}$ is closed and convex.
\end{corollary}

We now consider the simple, but non-trivial case that $\left( f_{i}\right)
_{i\in I}$ is a family of affine functions with a proper robust sum $f.$

\begin{example}
\label{Exam2}Let
\begin{equation*}
f_{i}=\left\langle a_{i}^{\ast },\cdot \right\rangle -t_{i},\ \left(
a_{i}^{\ast },\ t_{i}\right) \in X^{\ast }\times \mathbb{R},\ \forall i\in I,
\end{equation*}%
and suppose that there exist $\overline{x}\in X$ and $M\in \mathbb{R}$ such
that
\begin{equation}
\sum\limits_{i\in J}\left( \left\langle a_{i}^{\ast },\overline{x}%
\right\rangle -t_{i}\right) \leq M,\ \forall J\in \mathcal{F}\left( I\right)
.  \label{5.4}
\end{equation}%
\newline
For each $i\in I,$ we have $f_{i}^{\ast }=\delta _{a_{i}^{\ast }}+t_{i},$
where $\delta _{a_{i}^{\ast }}:X^{\ast }\longrightarrow \mathbb{R\cup }%
\left\{ +\infty \right\} $ represents the indicator function of $a_{i}^{\ast
},$ i.e., $\delta _{a_{i}^{\ast }}\left( x^{\ast }\right) =0,$ if $x^{\ast
}=a_{i}^{\ast },$ and $\delta _{a_{i}^{\ast }}\left( x^{\ast }\right)
=+\infty ,$ otherwise. \ Defining $A:\mathcal{F}\left( I\right)
\longrightarrow X^{\ast }$ such that $A\left( J\right) =\sum\nolimits_{i\in
J}a_{i}^{\ast },$ the function $\varphi $ writes%
\begin{equation*}
\varphi \left( x^{\ast }\right) :=\inf_{J\in A^{-1}\left( x^{\ast }\right)
}\dsum\limits_{i\in J}t_{i},\ \forall x^{\ast }\in X^{\ast }.
\end{equation*}%
\newline
The robust sum problem $(\mathrm{RP}_{x^{\ast }})$ has a zero duality gap
means that
\begin{equation}
\inf_{x\in X}\sup_{J\in \mathcal{F}\left( I\right) }\dsum\limits_{i\in
J}\left( \left\langle a_{i}^{\ast }-{x}^{\ast },x\right\rangle -t_{i}\right)
=\sup_{J\in A^{-1}\left( \overline{x}^{\ast }\right) }-\dsum\limits_{i\in
J}t_{i}.  \label{5.2}
\end{equation}%
\newline
We note that, given $\varepsilon \geq 0,$
\begin{equation*}
M^{\varepsilon }f_{i}(x^{\ast })=\left\{
\begin{array}{ll}
X,\ \ \ \ \ \ \ \ \  & \mathrm{if}\ x^{\ast }=a_{i}^{\ast }, \\
\emptyset , & \text{\textrm{else}}.%
\end{array}%
\right.
\end{equation*}%
\newline
Consequently, from (\ref{3}),
\begin{equation}
N^{\varepsilon }f(x^{\ast })=\bigcap_{\eta >0}\dbigcup\limits_{\QATOP{J\in
A^{-1}\left( \overline{x}^{\ast }\right) }{\hfill }}T_{f}^{\varepsilon +\eta
}\left( J\right) .  \label{5.3}
\end{equation}%
By Theorem \ref{Th1}, (\ref{5.2}) holds if and only if
\begin{equation*}
M^{\varepsilon }f(\overline{x}^{\ast })=\bigcap_{\eta >0}\dbigcup\limits_{%
\QATOP{J\in A^{-1}\left( \overline{x}^{\ast }\right) }{\hfill }%
}T_{f}^{\varepsilon +\eta }\left( J\right) ,\ \forall \varepsilon \geq 0.
\end{equation*}%
By (\ref{5.3}) one has $x^{\ast }\in \left( N^{\varepsilon }f\right)
^{-1}(x) $ if and only if for each $\eta >0$ there exists $J\in \mathcal{F}%
\left( I\right) $ such that $x^{\ast }\in A\left( J\right) $ and $J\in
S_{f}^{\varepsilon +\eta }\left( x\right) ,$ that means
\begin{equation*}
x^{\ast }\in \bigcap_{\eta >0}\dbigcup\limits_{\QATOP{J\in
S_{f}^{\varepsilon +\eta }\left( x\right) }{\hfill }}\sum\limits_{i\in
J}a_{i}^{\ast }.
\end{equation*}%
Consequently, from Theorem \ref{Th2},\ (\ref{5.2}) holds for each $x^{\ast
}\in X^{\ast }$\ if and only if
\begin{equation*}
\partial ^{\varepsilon }f\left( x\right) =\bigcap_{\eta >0}\dbigcup\limits_{%
\QATOP{J\in S_{f}^{\varepsilon +\eta }\left( x\right) }{\hfill }%
}\sum\limits_{i\in J}a_{i}^{\ast },\ \forall \left( x,\varepsilon \right)
\in X\times \mathbb{R}_{+}.
\end{equation*}%
Regarding the closedness criteria in Theorem \ref{Th5} and Corollary \ref%
{Corol1}, observe that%
\begin{equation*}
\begin{array}{ll}
\mathcal{A} & =\dbigcup\limits_{\QATOP{J\in \mathcal{F}\left( I\right) }{%
\hfill }}\sum\limits_{i\in J}\limfunc{epi}f_{i}^{\ast } \\
& =\dbigcup\limits_{J\in \mathcal{F}\left( I\right) }\left[ \left\{
\sum\limits_{i\in J}\left( a_{i}^{\ast },t_{i}\right) \right\} +\left\{
0_{X^{\ast }}\right\} \times \mathbb{R}_{+}\right]%
\end{array}%
\end{equation*}%
is the union of infinitely many vertical closed half-lines.

It is worth observing that in case all functions are linear (i.e., $t_{i}=0$
for all $i\in I$), $\mathcal{A}$ is closed (convex, respectively) if and
only if $\left\{ \sum\limits_{i\in J}a_{i}^{\ast }:J\in \mathcal{F}\left(
I\right) \right\} $ is closed (convex). When $I$ is countable (as in the
robust sums of linear functions in the third example of the introduction), $%
\left\{ \sum\limits_{i\in J}a_{i}^{\ast }:J\in \mathcal{F}\left( I\right)
\right\} $ is countable too, so that it cannot be convex. Finally, in the
simplest case that all functions are constants (i.e., $a_{i}^{\ast
}=0_{X^{\ast }}$ for all $i\in I$ and, according to (\ref{5.4}), $\theta
:=\sum\nolimits_{i\in I}^{R}-t_{i}\in \mathbb{R}$), we have $\{0_{X^{\ast
}}\}\times ]-\theta ,+\infty \lbrack \ \subset \mathcal{A}\subset
\{0_{X^{\ast }}\}\times \lbrack -\theta ,+\infty \lbrack $ and either $%
\mathcal{A}=\{0_{X^{\ast }}\}\times ]-\theta ,+\infty \lbrack $ or $\mathcal{%
A}=\{0_{X^{\ast }}\}\times \lbrack -\theta ,+\infty \lbrack $. So, $\mathcal{%
A}$ is convex. However, $\mathcal{A}$ is closed if and only if there exists $%
J\in \mathcal{F}(I)$ such that $\theta =\sum\limits_{i\in J}-t_{i}$.
\end{example}

\section{Duality for the infinite sum of non-negative convex functions and
related situations}

The case that the functions $f_{i},$ $i\in I,$ are non-negative presents
many specificities. For example, in such a case the robust sum coincides
with the infinite sum, i.e.,%
\begin{equation*}
f{\left( x\right) =}\sum\nolimits_{i\in I}^{R}{{f_{i}}\left( x\right) =}%
\lim\limits_{J\in \mathcal{F}\left( I\right) }\sum\limits_{i\in J}{{f_{i}}}%
\left( x\right) =\sum\limits_{i\in I}{{f_{i}}}\left( x\right) ,\forall x\in
X.
\end{equation*}
As in Section 2, the limit is taken respect to the directed set $\mathcal{F}%
\left( I\right) $ ordered by the inclusion relation. We have the next
important convexity properties.

\begin{lemma}
\label{Lemma8}Assume that $f_{i}\geq 0$ for each $i\in I.$ Then the set $%
\mathcal{A}=\dbigcup\limits_{\QATOP{J\in \mathcal{F}\left( I\right) }{\hfill
}}\sum\limits_{i\in J}\limfunc{epi}f_{i}^{\ast }$ and the function $\varphi $
are convex.
\end{lemma}

\textit{Proof.} Let $\left( x^{\ast },r\right) ,\left( y^{\ast },s\right)
\in \mathcal{A}$ and $t\in \left[ 0,1\right] .$ There exist $J,K\in \mathcal{%
F}\left( I\right) $ such that $\left( x^{\ast },r\right) \in
\sum\nolimits_{j\in J}\limfunc{epi}f_{j}^{\ast }$ and $\left( y^{\ast
},s\right) \in \sum\nolimits_{k\in K}\limfunc{epi}f_{k}^{\ast }.$

Let $l\in I.$ As $f_{l}\geq 0,$ we have $f_{l}^{\ast }\left( 0_{X^{\ast
}}\right) \leq 0,$ that is $\left( 0_{X^{\ast }},0\right) \in \limfunc{epi}%
f_{l}^{\ast }.$ Let $L:=J\cup K\in \mathcal{F}\left( I\right) .$ Since $%
\left( 0_{X^{\ast }},0\right) \in \limfunc{epi}f_{l}^{\ast }$ for all $l\in
L,$ $\left( x^{\ast },r\right) ,\left( y^{\ast },s\right) \in
\sum\nolimits_{l\in L}\limfunc{epi}f_{l}^{\ast } $ which is a convex subset
of $\mathcal{A}$. Thus, $\left( 1-t\right) \left( x^{\ast },r\right)
+t\left( y^{\ast },s\right) \in \mathcal{A}$ and $\mathcal{A}$ is convex.

The convexity of $\varphi $ is \ a consequence of (\ref{5.1}). In fact, we
have%
\begin{equation*}
\varphi \left( x^{\ast }\right) =\inf \left\{ r\in \mathbb{R}:\left( x^{\ast
},r\right) \in \mathcal{A}\right\} ,\forall x^{\ast }\in X^{\ast },
\end{equation*}%
which is a convex functions thanks to the convexity of $\mathcal{A}.\qquad $%
\qed

In what follows we assume that%
\begin{equation}
\left( f_{i}\right) _{i\in I}\subset \Gamma \left( X\right) ,\text{ }%
f=\sum\nolimits_{i\in I}^{R}{{f_{i}}}\text{ is proper, and }\mathcal{A}%
=\dbigcup\limits_{\QATOP{J\in \mathcal{F}\left( I\right) }{\hfill }%
}\sum\nolimits_{i\in J}\limfunc{epi}f_{i}^{\ast }\text{ is convex.}
\label{8}
\end{equation}

\begin{lemma}
\label{Lemma9}Assume that (\ref{8}) holds. Then $f^{\ast }=\overline{\varphi
}$ (the $w^{\ast }$-lsc hull of $\varphi $).
\end{lemma}

\textit{Proof.} By Lemma \ref{Lemma7}, we have $\varphi ^{\ast \ast
}=f^{\ast }.$ As shown in the proof of Lemma \ref{Lemma8}, $\varphi $ is
convex due to the convexity of $\mathcal{A}$. Since $f$ is proper, one has $%
\func{dom}\varphi ^{\ast }\neq \emptyset $ and, consequently, $\overline{%
\varphi }=\varphi ^{\ast \ast }=f^{\ast }.\qquad $\qed

\begin{lemma}
\label{Lemma10} Assume that (\ref{8}) holds. Then for any $x\in X$ and any $%
\varepsilon >0,$ we have%
\begin{equation*}
\partial ^{\varepsilon }f\left( x\right) =\overline{\Pi _{s}^{\varepsilon
}f\left( x\right) }.
\end{equation*}
\end{lemma}

\textit{Proof.} If $f\left( x\right) =+\infty $, then $\partial
^{\varepsilon }f\left( x\right) =\overline{\Pi _{s}^{\varepsilon }f\left(
x\right) }=\emptyset .$ Assume now $f\left( x\right) \in \mathbb{R}.$ By
Lemma \ref{Lemma9}, $f^{\ast }=\overline{\varphi }$ and it now follows from %
\eqref{eqsubdif} that
\begin{equation}  \label{eq7.2a}
\partial ^{\varepsilon }f\left( x\right) =\left[ \overline{\varphi }%
-\left\langle \cdot ,x\right\rangle +f\left( x\right) \leq \varepsilon %
\right] =\left[ \overline{\varphi -\left\langle \cdot ,x\right\rangle
+f\left( x\right) }\leq \varepsilon \right] .
\end{equation}%
As $\varphi ^{\ast }\left( x\right) =f\left( x\right)$ (by Lemma \ref{Lemma7}%
), we have $0 = f(x) - \varphi ^\ast (x) = \inf\limits_{X^{\ast }}\left\{
\varphi -\left\langle \cdot ,x\right\rangle +f\left( x\right)
\right\}~<~\varepsilon.$ By \cite[Lemma 1.1]{HMSV95} (applies to the
function $\varphi -\left\langle \cdot ,x\right\rangle +f\left( x\right))$ we
have $\left[ \overline{{\varphi }-\left\langle \cdot ,x\right\rangle
+f\left( x\right)} \leq \varepsilon \right] = \overline{\left[ \varphi
-\left\langle \cdot ,x\right\rangle +f\left( x\right) <\varepsilon \right] }%
. $ Taking \eqref{eq7.2a} into account, we have
\begin{equation*}
\partial ^{\varepsilon }f\left( x\right) = \overline{\left[ \varphi
-\left\langle \cdot ,x\right\rangle +f\left( x\right) <\varepsilon \right] }.
\end{equation*}
Now it is straightforward to check that $\left[ \varphi -\left\langle \cdot
,x\right\rangle +f\left( x\right) <\varepsilon \right] \subset \Pi
_{s}^{\varepsilon }f\left( x\right) $, and hence,
\begin{equation*}
\partial ^{\varepsilon }f\left( x\right) = \overline{\left[ \varphi
-\left\langle \cdot ,x\right\rangle +f\left( x\right) <\varepsilon \right] }
\ \subset\ \overline{ \Pi _{s}^{\varepsilon }f\left( x\right)}.
\end{equation*}%
It now follows from \eqref{eq5.4} and Lemma \ref{Lemma6},%
\begin{equation*}
\Pi _{s}^{\varepsilon }f\left( x\right)\ \subset \ \Pi ^{\varepsilon
}f\left( x\right)\ \subset\ \partial ^{\varepsilon }f\left( x\right)\
\subset \ \overline{\Pi _{s}^{\varepsilon }f\left( x\right) }.
\end{equation*}%
Since $\partial ^{\varepsilon }f\left( x\right) $ is $w^{\ast }$-closed, we
get $\partial ^{\varepsilon }f\left( x\right) =\overline{\Pi
_{s}^{\varepsilon }f\left( x\right) }.\qquad $\qed

\begin{theorem}[Stable zero duality gap under convexity]
\label{Th6}Assume that (\ref{8}) holds. The next statements are equivalent:%
\newline
$(i)$ The robust sum problem $(\mathrm{RP}_{x^*}) $ has stable zero duality
gap, \newline
$(ii)$ $\Pi ^{\varepsilon }f\left( x\right) =\overline{\Pi _{s}^{\varepsilon
}f\left( x\right) }, \ \forall x\in X, \ \forall \varepsilon >0$, \newline
$(iii)$ There exists $\overline{\varepsilon }>0$ such that%
\begin{equation*}
\Pi ^{\varepsilon }f\left( x\right) =\overline{\Pi _{s}^{\varepsilon
}f\left( x\right) }, \ \forall \left( x,\varepsilon \right) \in X\times %
\left] 0,\overline{\varepsilon }\right[ ,
\end{equation*}%
\newline
$(iv)$ There exists $\delta >0$ such that%
\begin{equation*}
\partial ^{\varepsilon }f\left( x\right) \subset \Pi _{s}^{\varepsilon
\delta }f\left( x\right) ,\forall \left( x,\varepsilon \right) \in X\times
\left] 0,+\infty \right[ .
\end{equation*}
\end{theorem}

\textit{Proof.} The equivalence of $(i),$ $(ii)$ and $(iii)$ follows from
Theorem \ref{Th2} and Lemma \ref{Lemma10}.

$\left[ (i)\Longrightarrow (iv)\right] $ By Theorem \ref{Th2} we have $%
\partial ^{\varepsilon }f\left( x\right) =\Pi ^{\varepsilon }f\left(
x\right) .$ Now%
\begin{equation*}
\Pi ^{\varepsilon }f\left( x\right) =\bigcap_{\eta >0}\Pi _{s}^{\varepsilon
+\eta }f\left( x\right) \subset \Pi _{s}^{2\varepsilon }f\left( x\right)
\end{equation*}%
and $(iv)$ holds with $\delta =2.$

$\left[ (iv)\Longrightarrow (i)\right] $ Assume that $(i)$ does not hold and
let $\delta >0.$ We will show that there exist $x\in X$ and $\varepsilon >0$
such that%
\begin{equation}  \label{eq7.a}
\overline{\Pi _{s}^{\varepsilon }f\left( x\right) }\nsubseteq \Pi
_{s}^{\varepsilon \delta }f\left( x\right) .
\end{equation}%
Since $(i)$ does not hold, there exist $x^{\ast }\in X^{\ast }$ and $%
\varepsilon >0$ such that $f^{\ast }\left( x^{\ast }\right) +\varepsilon
\delta <\varphi \left( x^{\ast }\right) .$ Pick $\overline{x}\in \partial
^{\varepsilon }f^{\ast }\left( x^{\ast }\right) $ (which is non-empty since $%
\varepsilon >0$). We have $x^{\ast }\in \partial ^{\varepsilon }f\left(
\overline{x}\right) =\overline{\Pi _{s}^{\varepsilon }f\left( \overline{x}%
\right) }.$ Assume that $x^{\ast }\in \Pi _{s}^{\varepsilon \delta }f\left(
\overline{x}\right) .$ Then, exist $\alpha \in \left[ 0,\varepsilon \delta %
\right] ,$ $J\in S_{f}^{\alpha }\left( \overline{x}\right) ,$ $\left(
\varepsilon _{i}\right) _{i\in J}\in \mathbb{R}_{+}^{J},$ and $x_{i}^{\ast
}\in \partial ^{\varepsilon _{i}}f_{i}\left( \overline{x}\right) $ for all $%
i\in J,$ such that $\alpha +\sum\nolimits_{i\in J}\varepsilon
_{i}=\varepsilon \delta $ and $\sum\nolimits_{i\in J}x_{i}^{\ast }=x^{\ast
}. $ Then%
\begin{equation*}
\begin{array}{ll}
\varphi \left( x^{\ast }\right) & \leq \sum\limits_{i\in J}{{f_{i}^{\ast }}}%
\left( x_{i}^{\ast }\right) \\
& \leq \sum\limits_{i\in J}\left( \left\langle x_{i}^{\ast },\overline{x}%
\right\rangle -{{f_{i}}}\left( \overline{x}\right) +\varepsilon _{i}\right)
\\
& =\left\langle x^{\ast },\overline{x}\right\rangle -\sum\limits_{i\in J}{{%
f_{i}}}\left( \overline{x}\right) +\varepsilon \delta -\alpha \\
& \leq \left\langle x^{\ast },\overline{x}\right\rangle -{f}\left( \overline{%
x}\right) +\alpha +\varepsilon \delta -\alpha \\
& \leq f^{\ast }\left( x^{\ast }\right) +\varepsilon \delta <\varphi \left(
x^{\ast }\right) ,%
\end{array}%
\end{equation*}%
which contradicts $f^{\ast }\left( x^{\ast }\right) +\varepsilon \delta
<\varphi \left( x^{\ast }\right)$. So $x^{\ast }\notin \Pi _{s}^{\varepsilon
\delta }f\left( \overline{x}\right)$, \eqref{eq7.a} is proved and the proof
is complete. \qquad \qed

\begin{corollary}
\label{Corol2}Assume that (\ref{8}) holds and $\Pi _{s}^{\varepsilon
}f\left( x\right) $ is $w^{\ast }$-closed for each $x\in X$ and $\varepsilon
>0.$ Then the robust sum problem $(\mathrm{RP}_{x^*}) $ has stable zero
duality gap.
\end{corollary}

\textit{Proof.} Under the assumption, it follows from Lemma \ref{Lemma10}
that $\partial ^{\varepsilon }f\left( x\right) = \overline{\Pi
_{s}^{\varepsilon }f\left( x\right) } = \Pi _{s}^{\varepsilon }f\left(
x\right)$, which means that statement $(iv)$ of Theorem \ref{Th6} holds with
$\delta =1.$ \qquad \qed

\textbf{Acknowledgements} This research was supported by the National
Foundation for Science \& Technology Development (NAFOSTED), Vietnam,
Project \textit{Some topics on systems with uncertainty and robust
optimization}, and by the Ministry of Economy and Competitiveness of Spain
and the European Regional Development Fund (ERDF) of the European
Commission, Project MTM2014-59179-C2-1-P.

\end{document}